\date{}
\documentclass [11pt]{article}
\usepackage{amsfonts,amsmath,amssymb,color,graphicx} 

\title{\vspace{-0.5in} On a conjecture of Erd\H{o}s and Simonovits: Even cycles.}
\author{
Peter Keevash\thanks{
School of Mathematical Sciences, Queen Mary University of London, UK, p.keevash@qmul.ac.uk.
Research supported in part by ERC grant 239696 and EPSRC grant EP/G056730/1.}
\and Benny Sudakov\thanks{
Department of Mathematics, University of California, Los Angeles, CA, USA, bsudakov@math.ucla.edu.
Research supported in part by NSF grant DMS-1101185, NSF CAREER award DMS-0812005 and by
USA-Israeli BSF grant.
}
\and Jacques Verstra\"ete\thanks{
Department of Mathematics, University of California, San Diego, CA, USA, jverstra@math.ucsd.edu.
Research supported in part by an Alfred P. Sloan Research Fellowship and NSF Grant DMS-0800704.}
}

\newtheorem{thm}{Theorem}[section] 
\newtheorem{prop}[thm]{Proposition}
\newtheorem{lem}[thm]{Lemma}

\newcommand{\mc}[1]{\mathcal{#1}}

\newcommand{\nib}[1]{\noindent {\bf #1}}

\newcommand{\sm}{\setminus}
\newcommand{\ov}{\overline}

\newcommand{\eps}{\varepsilon}
\newcommand{\sub}{\subseteq}
\newcommand{\brac}[1]{\left( #1 \right)}
\newcommand{\bsize}[1]{\left| #1 \right|}

\def\qed{\hfill $\Box$}

\begin{document}

\maketitle

\begin{abstract}
Let $\mc{F}$ be a family of graphs.
A graph is {\em $\mc{F}$-free} if it contains no copy of a graph in $\mc{F}$ as a subgraph.
A cornerstone of extremal graph theory is the study of the {\em Tur\'{a}n number} $\mbox{ex}(n,\mc{F})$,
the maximum number of edges in an $\mc{F}$-free graph on $n$ vertices.
Define the {\em Zarankiewicz number} $\mbox{z}(n,\mc{F})$ to be the maximum number of edges
in an $\mc{F}$-free {\em bipartite} graph on $n$ vertices.
Let $C_k$ denote a cycle of length $k$, and let $\mc{C}_k$ denote the set of
cycles $C_{\ell}$, where $3 \le \ell \leq k$ and $\ell$ and $k$ have the same parity.
Erd\H{o}s and Simonovits conjectured that for any family $\mc{F}$ consisting of bipartite graphs
there exists an odd integer $k$ such that $\mbox{ex}(n,\mc{F} \cup \mc{C}_k) \sim \mbox{z}(n,\mc{F})$.
They proved this when $\mc{F}=\{C_4\}$ by showing that $\mbox{ex}(n,\{C_4,C_5\}) \sim \mbox{z}(n,C_4)$.
In this paper, we extend this result by showing that if $\ell \in \{2,3,5\}$ and $k > 2\ell$ is odd, then
$\mbox{ex}(n,\mc{C}_{2\ell} \cup \{C_k\}) \sim \mbox{z}(n,\mc{C}_{2\ell})$. Furthermore, if $k > 2\ell + 2$ is odd, then
for infinitely many $n$ we show that the extremal $\mc{C}_{2\ell} \cup \{C_k\}$-free graphs are bipartite incidence graphs
of generalized polygons. We observe that this exact result does not hold for any odd $k < 2\ell$,
and furthermore the asymptotic result does not hold when $(\ell,k)$ is $(3,3)$, $(5,3)$ or $(5,5)$.
Our proofs make use of pseudorandomness properties of nearly extremal graphs that are of independent interest.
\end{abstract}

\section{Introduction}

Given a family $\mc{F}$ of graphs, a graph is {\em $\mc{F}$-free} if it contains no copy of a graph in $\mc{F}$ as a subgraph.
The {\em Tur\'an number} $\mbox{ex}(n,\mc{F})$ is the maximum number of edges in an $\mc{F}$-free graph on $n$ vertices.
When $\mc{F}=\{F\}$ consists of a single forbidden graph we denote the Tur\'an number by $\mbox{ex}(n,F)$.
A classical theorem of Tur\'an \cite{T} gives an exact result for $\mbox{ex}(n,K_t)$, where
$K_t$ is the complete graph on $t$ vertices: the unique largest $K_t$-free graph
on $n$ vertices is the complete $(t-1)$-partite graph with part sizes as equal
as possible. In general, Erd\"os, Stone and Simonovits \cite{ESi2,ES} showed that
$\mbox{ex}(n,\mc{F}) = (1-1/r)\binom{n}{2}+o(n^2)$, where $r = \min \{ \chi(F)-1: F \in \mc{F} \}$.
This determines the Tur\'an number asymptotically when $\mc{F}$ consists of non-bipartite graphs.
However, much less is known concerning the Tur\'an numbers of bipartite graphs.
There is no bipartite graph $F$ containing a cycle such that $\mbox{ex}(n,F)$ is known exactly for all $n$.
In fact, even the order of magnitude of $\mbox{ex}(n,F)$ is not known for quite simple bipartite graphs,
such as the complete bipartite graph with four vertices in each part,
the cycle of length eight, and the three-dimensional cube graph.

In this paper, we study the effect of forbidding short odd cycles on bipartite extremal problems, in particular,
the extremal problem for even cycles.  Let $C_k$ denote a cycle of length $k$ and let $\mc{C}_k$ denote the set of
cycles $C_{\ell}$ where $3 \le \ell \leq k$ and $\ell$ and $k$ have the same parity.
Bondy and Simonovits~\cite{BS} showed that $\mbox{ex}(n,C_{2\ell}) = O(n^{1 + 1/\ell})$, which was further improved by Lam and Verstra\"ete~\cite{LV} to
\begin{equation}\label{eq:lv}
\mbox{ex}(n,\mc{C}_{2\ell}) \leq \frac{1}{2}n^{1 + 1/\ell} + O(n).
\end{equation}
It is a notoriously difficult problem to determine even the order of magnitude of $\mbox{ex}(n,\mc{C}_{2\ell})$
or $\mbox{z}(n,\mc{C}_{2\ell})$ for $\ell \not \in \{2,3,5\}$. For $\ell \in \{2,3,5\}$, the existence of
structures from projective geometry called {\em generalized polygons} show $\mbox{ex}(n,\mc{C}_{2\ell}) = \Theta(n^{1 + 1/\ell})$.
Erd\H{o}s and Simonovits~\cite{ESi1} conjectured more generally that  $\mbox{ex}(n,\mc{C}_{2\ell}) = \Theta(n^{1 + 1/\ell})$ for
all $\ell \geq 2$, and this conjecture remains open. For large $\ell$, the densest known $\mc{C}_{2\ell}$-free graphs on $n$ vertices
are the recent constructions of Ramanujan graphs based on octonions, due to Dahan and Tillich~\cite{DT}, superseding earlier constructions of Margulis~\cite{M}, Lubotzky, Phillips, and Sarnak~\cite{LPS}, and Lazebnik, Ustimenko and Woldar~\cite{LUW1}.
The constructions have $n^{1 + \theta}$ edges where $\theta \sim \frac{6}{7\ell}$ as $\ell \rightarrow \infty$. For $\ell = 4$, the current
best bounds are $c_1 n^{6/5} \leq \mbox{ex}(n,\mc{C}_{2\ell}) \leq c_2 n^{5/4}$ for some constants $c_1,c_2 > 0$ -- the lower bound is from
the constructions of Benson \cite{Ben} and Singleton \cite{S}, whereas the upper bound follows from (\ref{eq:lv}) with $\ell = 4$. In this paper, we study the relationship between extremal $\mc{C}_{2\ell} \cup \{C_k\}$-free graphs with $k$ odd, and extremal $\mc{C}_{2\ell}$-free bipartite graphs.

\subsection{Main Results}

Define the {\em Zarankiewicz number} $\mbox{z}(n,\mc{F})$ to be the maximum number of edges
in an $\mc{F}$-free bipartite graph on $n$ vertices.
Erd\H{o}s and Simonovits~\cite[Conjecture 3]{ESi1} conjectured that for any family $\mc{F}$ consisting of bipartite graphs
there exists an odd integer $k$ such that $\mbox{ex}(n,\mc{F} \cup \mc{C}_k) \sim \mbox{z}(n,\mc{F})$.
They proved this when $\mc{F} =\{C_4\}$ by showing that $\mbox{ex}(n,\{C_4,C_5\}) \sim \mbox{z}(n,\{C_4\})$.
In this paper, we extend this by proving the following theorem which verifies the Erd\H{o}s-Simonovits conjecture when
$\mc{F} = \mc{C}_{2\ell}$ for $\ell \in \{2,3,5\}$.

\begin{thm} \label{stability}
Suppose $\ell \in \{2,3,5\}$ and $k > 2\ell$ is odd.
Let $G$ be a $\mc{C}_{2\ell} \cup \{C_k\}$-free graph on $n$ vertices with average degree $d=\Theta(n^{1/\ell})$.
Then $G$ has a bipartite subgraph $H$ with at least $d^{\ell+1}-o(n^{1+1/\ell})$ edges.
If also $d \ge (1+o(1))(n/2)^{1/\ell}$ then $d \sim (n/2)^{1/\ell}$ and $e(H) \sim e(G) \sim (n/2)^{1 + 1/\ell}$.
In particular
\[ \mbox{ex}(n,\mc{C}_{2\ell} \cup \{C_k\}) \sim \mbox{z}(n,\mc{C}_{2\ell}).\]
\end{thm}

Note that the statement in Theorem \ref{stability} is stronger
than that of the conjecture in two ways. Firstly, we replace $\mc{C}_k$ by $C_k$,
i.e.\ we forbid a single odd cycle rather than all short odd cycles.
Secondly, we obtain a {\em stability theorem}: any graph which is close to extremal
for $\mc{C}_{2\ell} \cup \{C_k\}$ must be close to bipartite. This is the first ingredient
in applying the {\em stability method}, in which one first obtains approximate structure,
and then eliminates any imperfections to obtain exact structure.
Using this theorem  we can prove the following exact result. The reader unfamiliar with generalized
$(\ell + 1)$-gons can find a brief description in Section \ref{sec:zarankiewicz}: these can be viewed for infinitely many prime
powers $q$ as $(q + 1)$-regular $\mc{C}_{2\ell}$-free extremal bipartite graphs with $q^{\ell} + q^{\ell - 1} + \dots + 1$ vertices in each part.

\begin{thm}\label{even}
Let $\ell \in \{2,3,5\}$, $k > 2\ell + 2$ be odd, $n>n_k$ be sufficiently large,
and define $q \in \mathbb R^+$ by $n=2(q^\ell + q^{\ell - 1} + \dots + 1)$.
Then any $\mc{C}_{2\ell}$-free graph $G$ on $n$ vertices with at least $\frac{1}{2}(q + 1)n$ edges
contains a cycle of length $k$, unless $q$ is an integer and $G$ is the bipartite incidence graph
of a generalized $(\ell + 1)$-gon of order $q$. Furthermore, if $q$ is an integer, $n > n_k$ and there
is a generalized $(\ell + 1)$-gon of order $q$, then
\[ \mbox{ex}(n,\mc{C}_{2\ell} \cup \{C_k\}) = \mbox{z}(n,\mc{C}_{2\ell}).\]
The same statement holds when $\ell = 2$ and $k = 5$.
\end{thm}

\subsection{Chromatic number}

A classical result of Andr\'asfai, Erd\H{o}s and S\'{o}s~\cite{AES} states that
a triangle-free $n$-vertex graph with minimum degree more than $2n/5$ is 2-colorable. Generalizations of this
theorem have been studied extensively by researchers, for example see~\cite{Al, ASu, GL, LT} and their references.
Here we address a similar question
for a $\mathcal{C}_{2\ell} \cup \{C_k\}$-free graph when $k \geq 4\ell + 1$ is odd. We use $\chi(G)$
to denote the chromatic number of $G$.

\begin{thm}\label{chromatic}
Let $\ell \geq 2$ be an integer, and let $k \geq 4\ell + 1$ be an odd integer, let $c$ be a positive real number, and let $G$ be a $\mathcal{C}_{2\ell} \cup \{C_k\}$-free graph on $n$ vertices with minimum degree at least $cn^{1/\ell}$. Then $\chi(G) < (4k)^{\ell + 1}/c^{\ell}$.
\end{thm}

Forbidding short odd cycles in a graph generally has little effect on the chromatic number if the graph is too sparse.
A well-known construction of Erd\H{o}s in random graphs shows that there are graphs of arbitrarily large girth and chromatic number. Theorem \ref{chromatic} in contrast shows that the chromatic number becomes bounded for very dense $\mathcal{C}_{2\ell}$-free graphs with a forbidden long odd cycle.

\subsection{Short odd cycles}

We start by observing that the second statement of Theorem \ref{even}
does not hold whenever $5 \leq k < 2\ell$ and $k$ is odd,
as in this case we have
\[ \mbox{ex}(n,\mc{C}_{2\ell} \cup \{C_k\}) \geq \mbox{z}(n,\mc{C}_{2\ell}) + 1.\]
To see this, consider an extremal $\mc{C}_{2\ell}$-free bipartite graph $H$ on $n$ vertices.
Let $G$ be obtained from $H$ by adding an edge joining a pair of vertices $\{x,y\}$
at distance two in one part of $H$. We claim that $G$ has no $C_k$ for odd $k$ with $5 \leq k < 2\ell$.
For such a cycle would have to contain the edge $\{x,y\}$, so we would have a path $P$ of length $k-1$ in $H$ from $x$ to $y$.
Adding the edges $xz$ and $yz$, where $z$ is a common neighbor of $x$ and $y$, we obtain a closed walk of length $k+1$. Furthermore, this walk is not acyclic, so it must contain an even cycle of length at most $2\ell$ in $H$,
contradicting the fact that $H$ is $\mc{C}_{2\ell}$-free.
It would be interesting to see if Theorem \ref{even} can be extended to the last remaining case, namely $k = 2\ell + 1$.

The following proposition determines an upper bound for $\mbox{z}(n,\mc{C}_{2\ell})$.
Its proof will follow easily from the counting arguments in Section \ref{sec:zarankiewicz}.

\begin{prop} \label{evenbound}
Suppose $n \in \mathbb N$,
and let $q \in \mathbb R^+$ be defined by $n = 2(q^{\ell} + q^{\ell-1} +  \cdots + 1)$.
Then $\mbox{z}(n,\mc{C}_{2\ell}) \leq \frac{1}{2}(q + 1)n$.
\end{prop}

Now we can describe a much stronger discrepancy between Tur\'an and Zarankiewicz numbers when one
forbids a short odd cycle. Consider the {\em polarity graphs} constructed by Lazebnik, Ustimenko and
Woldar~\cite{LUW2}. These are $\mc{C}_{2\ell}$-free graphs on $n$ vertices with $(1/2+o(1))n^{1 + 1/\ell}$ edges,
such that for $\ell = 3$ they have no triangles and no even cycles of length at most six,
and for $\ell = 5$ they have no triangles and no cycles of length five and no even cycles of length at most ten.
Thus when $(\ell,k)$ is $(3,3)$, $(5,3)$ or $(5,5)$ we have
\[ \liminf_{n \rightarrow \infty} \frac{\mbox{ex}(n,\mc{C}_{2\ell} \cup \{C_k\})}{\mbox{z}(n,\mc{C}_{2\ell})}
\geq \liminf_{n \rightarrow \infty}  \frac{(n^{1 + 1/\ell})/2}{(n/2)^{1 + 1/\ell}} = 2^{1/\ell},\]
so in these cases one does not even have an asymptotic result similar to Theorem \ref{stability}.
For the case $\ell = 2$ and $k = 3$ we do not know of such a strong discrepancy.
Parsons \cite{P} constructed $\{C_3,C_4\}$-free graphs showing that
\[ \mbox{ex}(n,\{C_3,C_4\}) > \mbox{z}(n,C_4) + \frac{7}{32}n - O(\sqrt{n})\]
when $n = {q \choose 2}$ and $q = 1$ mod $4$ is prime.
On the other hand, Erd\H{o}s~\cite{ErdosA,ErdosB} suggested that in the case of $\{C_3,C_4\}$-free graphs there should not be a
stronger discrepancy, and conjectured that $\mbox{ex}(n,\{C_3,C_4\}) \sim \mbox{z}(n,C_4)$ -- this conjecture remains open.
One may also ask in general whether or not
$\mbox{ex}(n,\mc{C}_{2\ell} \cup \{C_{2\ell -1}\}) \sim \mbox{z}(n,\mc{C}_{2\ell})$ for $\ell \ge 2$.

\subsection{Organization}

This paper is organized as follows.
In the next section we illustrate our ideas by sketching the proofs of Theorems \ref{stability} and \ref{even} in the case of quadrilateral-free graphs.
Section 3 contains the essential results on counting walks and paths in graphs. These are used throughout the paper,
in particular to deduce some bounds on Tur\'{a}n numbers and Zarankiewicz numbers for cycles in the same section.
In Section 4 we show that nearly extremal graphs contain large subgraphs which are almost regular.
We prove a pseudorandomness property for nearly extremal graphs in Section 5.
Section 6 contains the proof of Theorem \ref{stability} and Section 7 the proof of Theorem \ref{even}. Section 8 contains
the short proof of Theorem \ref{chromatic}.
We make some concluding remarks in the final section.

\subsection{Notation}

We write $e(G)$ for the number of edges in a graph $G$.
In a graph $G$, let $N_r(v)$ denote the set of vertices at distance exactly $r$ from $v$, and let $d_r(v) = |N_r(v)|$.
For $r = 1$, we omit the subscript $r$, so that $d(v)$ is the degree of $v$ and $N(v)$ is the neighborhood of $v$.
We write $d_B(v)$ for the number of neighbors of a vertex $v$ in a set $B$.
We write $G[S]$ for the subgraph of $G$ induced by a set $S \sub V(G)$
and $e(S)$ for the number of edges in $G[S]$.
Given two sets $S,T \sub V(G)$, not necessarily disjoint, we write $e(S,T)$ for the number of ordered pairs
$(s,t)$ with $s \in S$, $t \in T$ and $st \in E(G)$. For example $e(S,S)=2e(S)=2|E(S)|$.
In addition to the {\em Tur\'an number} $\mbox{ex}(n,\mc{F})$ and {\em Zarankiewicz number} $\mbox{z}(n,\mc{F})$ defined above,
we use $\mbox{z}(a,b,\mc{F})$ for the maximum number of edges in an $\mc{F}$-free bipartite graph
that has $a$ vertices in one part and $b$ vertices in the other part.
We let $\mathbb R^+$ denote the positive reals and $\mathbb N$ the positive integers.
Our asymptotic notation assumes that $n \to \infty$; we write $f(n) \sim g(n)$ if
$\lim_{n \rightarrow \infty} f(n)/g(n) = 1$, and $f(n) = o(g(n))$ if $\lim_{n \rightarrow \infty} f(n)/g(n) = 0$.

\section{Sketch proof for $4$-cycles}\label{sec:c4}

In this section we outline the proofs of Theorems \ref{stability} and \ref{even}
in the case $\ell=2$, i.e.\ $4$-cycles. This introduces the main ideas of our
approach, without some technicalities that arise in the other cases.

\subsection{Stability}

We start with Theorem \ref{stability} for $\ell=2$.
The idea is that we can take $H$ to be the bipartite subgraph of $G$, containing all the edges between
$N_2(v)$ and $N_3(v)$ for a suitable vertex $v$.
Suppose that $G$ is a graph on $n$ vertices with average degree $d$
and $G$ does not contain $C_4$ or $C_k$ for some odd $k \ge 5$. The proof proceeds by the following steps.

\nib{1. Controlling the maximum degree.}
Lemma \ref{c4maxdeg}(i) will show that we can delete at most $2n$ edges from $G$
to obtain a subgraph $G'$ with maximum degree $\Delta \le 2\sqrt{n}$.

\nib{2. Enumeration of walks and paths.}
The Blakley-Roy inequality (Proposition \ref{walks}) gives a lower bound on the number of walks.
It implies that we can find a vertex $v$ that is the start
of at least $d'{}^3$ walks of length $3$, where $d' \ge d-4$
is the average degree of $G'$. Also, the bound on the maximum degree
implies that all but $O(n)$ of these walks are paths, i.e.\
we have $d^3-O(n)$ paths of length $3$.

\nib{3. Finding odd cycles.} We shall show in Section \ref{sec:findodd} that if $G[N_2(v)]$ has average degree at least $\max(6,2k-8)$, then $G$ contains $C_k$, so we conclude that $G[N_2(v)]$ has average degree less than $\max(6,2k - 8)$.
This implies that the number of paths of length 3 which start at $v$ and end in $N_2(v)$ is $O(n)$. Since $G$ is $C_4$-free, we also have no path of length 3 starting in $v$ and ending in $N_1(v)$. Therefore, all but $O(n)$ of the paths found in step 2 go from $v$ to $N_3(v)$.

\nib{4. Conclusion.}
Since $G$ is $C_4$-free, each edge between $N_2(v)$ and $N_3(v)$ is contained in at most
one path of length $3$ from $v$ to $N_3(v)$. Thus we obtain $d^3-O(n)$ edges between $N_2(v)$ and $N_3(v)$,
and these constitute the bipartite subgraph $H$ needed for the first statement in Theorem \ref{stability}.
For the second statement, note that Proposition \ref{evenbound} implies $z(n,C_4) \le (1+o(1))(n/2)^{3/2}$.
Since $H$ is a bipartite $C_4$-free graph we must have $d \le (1+o(1))(n/2)^{1/2}$.
If also $d \ge (1+o(1))(n/2)^{1/2}$ then we have $d \sim (n/2)^{1/2}$
and $e(H) \sim e(G) \sim (n/2)^{3/2}$, which proves Theorem \ref{stability} for $\ell=2$, apart from the
case $k = 5$ which is proved by slightly refining the above arguments.

\subsection{Exact result}\label{sketch-exact}

Now we sketch Theorem \ref{even} for $\ell=2$.
Suppose that $n$ is large, $G$ is a graph on $n$ vertices with $e(G) \ge (q+1)n/2$,
where $q \in \mathbb R^+$ is defined by $n=2(q^2+q+1)$,
and $G$ does not contain $C_4$ or $C_k$ for some odd $k \ge 7$.
The proof proceeds by the following steps.

\nib{1. Pseudorandomness.} Not much is known about the structure of nearly extremal
graphs for bipartite Tur\'an problems. Here we obtain a result in this direction.
It says that the number of edges between any two large sets in nearly extremal graph is close
to what one would expect in a random graph with the same edge density: more precisely
we shall show that if $G$ is a $C_4$-free bipartite graph on $n$ vertices with
parts $X$ and $Y$ and average degree $d \sim (n/2)^{1/2}$, then for any $S \sub X$ and $T \sub Y$ we have
$e(S,T) = \frac{2d}{n}|S||T| + o(n^{3/2})$. This result is a special case of Theorem \ref{bippr}.

\nib{2. Controlling the minimum degree.}
We reduce the proof to the case when the minimum
degree satisfies $\delta(G)> q/4$. This uses a vertex
deletion argument that is quite standard in extremal graph theory.
We consider a sequence of graphs $G=G_n,G_{n-1},\cdots,G_t$ for some $0 \le t \le n$,
where $G_{i-1}$ is obtained from $G_i$ by deleting a vertex of degree at most $q/4$, while possible.
Some calculations show that this process must terminate with $t>n/2$.
Now suppose that we know Theorem \ref{even} holds
under the additional assumption $\delta(G)> q/4$. Applying this to $G_t$
gives $e(G_t) \le (r+1)t/2$, where $r$ is defined by $t=2(r^2+r+1)$.
Furthermore, we have $e(G) \le e(G_t) + (n-t)q/4$,
and calculations show that this is less than $(q+1)n/2$, unless $t=n$.
Since $e(G) \ge (q+1)n/2$, we must have $t=n$, so $\delta(G)>q/4$,
and we are justified in assuming this when proving Theorem \ref{even}.

\begin{figure}
\begin{center}
\includegraphics[width=3in]{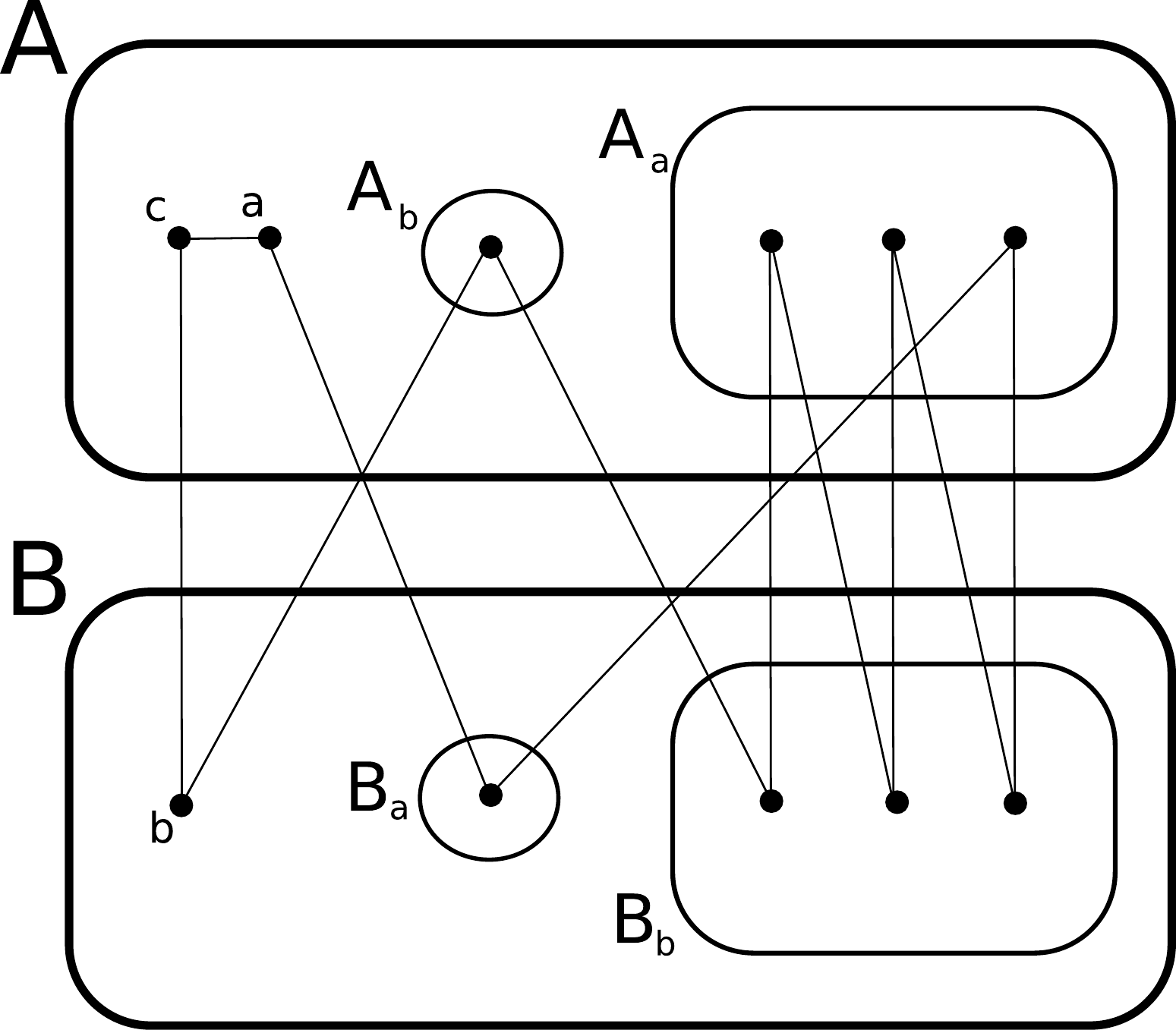}
\caption{Constructing a cycle of any specified odd length $k \ge 7$}
\end{center}
\end{figure}

\nib{3. Refining the approximate structure.}
Now we have $e(G) \ge (q+1)n/2$ and can assume that $\delta(G)>q/4$.
Let $H$ be a bipartite subgraph of $G$ with maximum size.
Theorem \ref{stability} implies that $e(H) \sim e(G) \sim (n/2)^{3/2}$.
Furthermore, maximality of $H$ implies that $\delta(H)>q/8$. Actually,
we only need the fact that $\delta(H)>cn^{1/2}$, where $c>0$ is independent of $n$.
Now we will show that $G=H$. Label the parts of $H$ as $A$ and $B$.
Suppose for a contradiction that $G[A]$ contains an edge $ac$.
Let $b$ be a neighbor of $c$ in $B$. Now we `explore' the graph
until we reach two sets of linear size where we can apply pseudorandomness (see Figure 1). By the minimum degree
assumption, we can take $|A_b|,
|B_a| \sim cn^{1/2}$ such that all vertices in $A_b$ are neighbors of $b$ and all vertices in $B_a$ are neighbors of $a$.
Next, taking all the neighbors of vertices in $A_b$ we get a set $B_b$ with $|B_b|= \Theta(n)$ and by taking neighbors of
$B_a$ we get a set $A_a$ with  $|A_a| = \Theta(n)$.
Then step 1 (Pseudorandomness) gives
\[ e(A_a,B_b) = \frac{2d}{n}|A_a||B_b| + o(n^{3/2}) = \Theta(n^{3/2}).\]
In particular, we can find a path of length $k-6$ using only edges between $A_a$ and $B_b$.
By construction this can be completed to a cycle of length $k$ in $G$, which is a contradiction.
We deduce that $A$ is an independent set in $G$. Similarly $B$ is independent, so $G=H$ is bipartite.
The characterization of equality now follows from a result on Zarankiewicz numbers for even cycles -- see
Proposition \ref{even=} -- so this proves Theorem \ref{even} for quadrilaterals.

\section{Counting walks}\label{sec:walks}

The basis of estimates for both Tur\'{a}n numbers and Zarankiewicz numbers for even cycles
is counting various types of walks in graphs. These counts are also used in Step 2 of the stability
result (Enumeration of walks and paths) and Step 1 of the exact result (Pseudorandomness).
A {\em walk of length $t$} in a graph $G$ is a sequence $v_0 e_0 v_1 e_1 \dots v_{t-1} e_{t-1} v_t$
such that $v_i \in V(G)$, $e_i \in E(G)$ and $e_i = \{v_i,v_{i + 1}\}$ for $0 \leq i < t$.
Note that edges may be repeated in a walk, even on consecutive steps.

\subsection{The Blakley-Roy Inequality}

Let $w_k(G)$ denote the number of walks of length $k$ in a graph $G$ divided by the number of vertices in $G$
-- this is the average number of walks of length $k$ starting at a vertex.
If $G$ is an $d$-regular graph on $n$ vertices, then clearly $w_k(G) = d^k$.
Blakley and Roy~\cite{BR} proved a matrix version of H\"{o}lder's Inequality,
which shows (as a special case) that any graph of average degree $d$
has at least as many walks of a given length as an $d$-regular graph on the same number of vertices:

\begin{prop}\label{walks}
Suppose $G$ is a graph of average degree $d$. Then $w_k(G) \geq d^k$.
\end{prop}

There are many proofs of this inequality in the literature;
the original proof of Blakley and Roy uses eigenvalues.
We now briefly discuss the tight connection between walks and eigenvalues.

\subsection{Walks and eigenvalues}

Let $G$ be a graph on $n$ vertices with adjacency matrix $A$,
and suppose $A$ has orthonormal eigenvectors $x_1,x_2,\dots,x_n$
and eigenvalues $\lambda_1 \ge \lambda_2 \ge \dots \ge \lambda_n$.
Then we may write
\[ w_k(G) = e^t A^k e \quad \mbox{ where } \quad e = n^{-1/2}(1,1,\dots,1).\]
In this notation, the Blakley-Roy inequality (Proposition \ref{walks}) mentioned earlier
can be stated as $w_k \ge w_1^k$. For future reference we note
the following proposition, attributed to Chris Godsil in \cite{ESi1}.

\begin{prop}\label{godsil}
Suppose $r,s \in \mathbb N$, where $r$ is even and $r \geq s$. Then $w_r^{1/r} \ge w_s^{1/s}$.
\end{prop}

\nib{Proof.}
Write $e = \sum c_i x_i$. Then $\sum c_i^2=e\cdot e = 1$ and $w_r = e^t A^r e =\sum c_i^2 \lambda_i^r$.
Jensen's inequality applied to $f(t)=t^{r/s}$ gives
$w_r = \sum c_i^2 \lambda_i^r = \sum c_i^2 |\lambda_i^s|^{r/s}
\ge \left( \sum c_i^2 |\lambda_i^s| \right)^{r/s} \ge w_s^{r/s}$.
\qed

\subsection{Closed walks and trace}

We can also use spectral theory to count closed walks.
Let $w^{\circ}_k$ be the number of closed walks of length $k$ in a graph $G$ divided by the number of vertices in $G$
-- this is the average number of closed walks of length $k$ starting at a vertex.

\begin{prop}\label{trace}
$w^{\circ}_k := \frac{1}{n}\mbox{Tr}(A^k) = \frac{1}{n}\sum_{i=1}^n \lambda_i^k$.
\end{prop}

This gives rise to a standard method for establishing a spectral gap in a $d$-regular graph:
if $w^{\circ}_k$ is close to $\lambda_1^k = d^k$, then all $|\lambda_i|$ for $i > 1$ must be small.
We have the following bound on $w^{\circ}_{2\ell + 2}$ in dense $\mc{C}_{2\ell}$-free
bipartite graphs with roughly equal part sizes.

\begin{lem}\label{closedwalks}
Suppose $G$ is a bipartite graph on $n$ vertices with part sizes $n/2 + o(n)$,
maximum degree $\Delta$ and girth at least $2\ell + 2$. Then
\[ w^{\circ}_{2\ell + 2}(G) < (1/2+o(1))n\Delta^2 + (4\Delta)^{\ell + 1}.\]
\end{lem}

\nib{Proof.}
Consider any vertex $v$ and a closed walk $W$ of length $2\ell + 2$ from $v$.
If $W$ is not a cycle then the girth assumption implies that its underlying graph is acyclic.
The number of such walks $W$ is therefore at most the number of closed walks of length $2\ell + 2$
from the root in the complete $\Delta$-ary tree. A crude upper bound on this number is
$\binom{2\ell+2}{\ell+1}\Delta^{\ell+1}<(4\Delta)^{\ell + 1}$, as may be seen by
choosing the $\ell+1$ times when the walk moves towards the root
and multiplying $\Delta$ choices for the $\ell+1$ times when the walk moves away from the root.
(The exact formula is $\frac{1}{\ell + 2}{2\ell + 2 \choose \ell + 1} \Delta^{\ell + 1}$
but we do not need this.) In the case when $W$ is a cycle, we estimate the possibilities
by considering the neighbors $a$ and $b$ of $v$ on $W$, and the opposite vertex $u$
at distance $\ell+1$ from $v$ on $W$. Note that $W$ is uniquely determined by $a$, $b$ and $u$,
as the girth assumption implies that there is at most one path of length $\ell$ between
any specified pair of vertices. We can choose each of $a$ and $b$ in at most $\Delta$ ways,
and $u$ in at most $n/2 + o(n)$ ways (using the assumption on the part sizes). Thus
the number of possibilities is at most $(1/2+o(1))n\Delta^2$, which also takes into
account the orientation of the cycle. Combining the two cases,
the number  of closed walks of length $2\ell + 2$ from $v$ is at most
$(1/2+o(1))n\Delta^2 + (4\Delta)^{\ell + 1}$. Since $v$ was arbitrary,
this bound also holds for the average $w^{\circ}_{2\ell + 2}(G)$.
\qed

\subsection{Non-returning walks}

A {\em non-returning walk of length $k$} in a graph $G$ is a walk $v_0 e_0 v_1 e_1 \dots v_{k-1} e_{k-1} v_k$
of length $k$ such that $e_i \neq e_{i + 1}$ for $0 \leq i < k$.
Let $\nu_k(G)$ denote the number of non-returning walks of length $k$ in a graph $G$ divided by the number of vertices in $G$.
If $G$ is an $r$-regular graph on $n$ vertices then clearly $\nu_k(G) \geq r(r - 1)^{k - 1}$.
Alon, Hoory and Linial~\cite{AHL} gave an analogue of Proposition \ref{walks} for non-returning walks.

\begin{prop}\label{nwalks}
Suppose $G$ is a graph of average degree $r \geq 2$. Then $\nu_k(G) \geq r(r - 1)^{k - 1}$. If $r \in \mathbb N$ then
equality holds if and only if $G$ is $r$-regular.
\end{prop}

In fact, they proved a slightly stronger form that the result given, which is sensitive to variations in the degrees of the vertices in the graph. Sidorenko~\cite{Sido} gave a bipartite analogue of Proposition \ref{walks}. Here we need a bipartite analogue of Proposition \ref{nwalks}, which was proved by Hoory \cite{H}. We will only need the result for walks of odd length.

\begin{prop}\label{bipnoreturn}
Suppose $G$ is a bipartite graph with parts $A$ and $B$ and average degree $d$.
Let $\alpha$ be the average degree of vertices in $A$
and $\beta$ the average degree of vertices in $B$.
Then for any $t \in \mathbb N$ we have
\[ \nu_{2t+1} \ge d \prod_{v \in A \cup B} (d(v) - 1)^{td(v)/e(G)} \geq d (\alpha - 1)^t(\beta - 1)^t.\]
If $\alpha,\beta \in \mathbb N$ then equality holds if and only if every vertex of $A$ has degree $\alpha$ and every
vertex of $B$ has degree $\beta$.
\end{prop}

\subsection{Paths}

A {\em path} is a walk which has no repeated vertices or edges. Let $p_\ell(G)$ denote the
number of paths of length $\ell$ in a graph $G$ divided by the number of vertices in $G$.

\begin{lem}\label{paths}
Suppose $G$ is a graph with maximum degree $\Delta$ and average degree $d$, and $\ell \in \mathbb N$. Then
\[ p_{\ell}(G) \geq d^{\ell} - \ell^2 \Delta^{\ell - 1}.\]
\end{lem}

\nib{Proof.}
By the Blakley-Roy inequality there are at least $nd^{\ell}$ walks of length $\ell$ in $G$ if $G$ has $n$ vertices.
Let $(v_1,v_2,\dots,v_{\ell})$ be a walk of length $\ell$ which is not a path.
Then $v_i = v_j$ for some distinct $i,j \in \{1,2,\dots,\ell\}$.
Fixing $i < j$, there are at most $n\Delta^{j-1}$ choices for the part of the walk up to $v_{j-1}$, then the next
step is determined since $v_j = v_i$, and then there are at most $\Delta^{\ell - j}$ choices for the remaining steps.
There are at most $\ell^2$ choices for $i$ and $j$, and therefore the number of walks of length $\ell$
which are not paths is at most $n\ell^2 \Delta^{\ell - 1}$. Dividing by $n$, we obtain the lemma.
\qed

\subsection{Finding odd cycles}\label{sec:findodd}

We use the following lemma which is implicit in \cite{V}.

\begin{lem}\label{find-odd}
Suppose $G$ is a graph, $v$ is a vertex of $G$ and $G[N_r(v)]$ has
average degree at least $2s - 4$, for some $r \ge 1$ and odd $s \ge 5$.
Then $G$ contains cycles $C_{2m + 1},C_{2m + 3},\dots,C_{2m + s}$ for some $1 \le m \le r$.
In particular, if $G$ is a graph not containing a cycle $C_k$, for some fixed odd $k>2r$, then
$G[N_r(v)]$ has average degree less than $\max(6,2k-8)$.
\end{lem}

\nib{Proof.}
To give a self-contained proof of this lemma would require the duplication of large parts of \cite{V},
so instead we just sketch the argument,
and refer the reader to \cite{V} for the omitted details.
Following the proof of \cite[Lemma 3]{V}, we note that $G[N_r(v)]$ has a subgraph with
minimum degree at least $s-1$, and so contains a subgraph $H$ which is a cycle of length at least $s$
with at least one chord. Next, as in the proof of \cite[Theorem 1]{V}, we consider a minimal subtree $T$ of
the breadth first search tree rooted at $v$ that contains $V(H)$. By minimality $T$ branches at its root.
We let $A$ be the set of vertices of $H$ belonging to one branch of $T$ and let $B=V(H) \sm A$.
Next we apply \cite[Lemma 2]{V}, which tells us that $H$ contains paths from $A$ to $B$ of every
length $\ell \le |V(H)|$, unless $H$ is bipartite with bipartition $(A,B)$. Even if $H$ is bipartite
with bipartition $(A,B)$ we still have paths from $A$ to $B$ of every odd length $\ell \le |V(H)|$.
These may be completed to cycles via paths in $T$, which all have the same length $2m$,
where $m$ is the distance from the root of $T$ to $N_r(v)$. The second statement of
the lemma is proved by taking $s = \max(5,k - 2)$. \qed

\subsection{Zarankiewicz numbers of even cycles}\label{sec:zarankiewicz}

In this section, we describe extremal bipartite $\mc{C}_{2\ell}$-free graphs and generalized polygons,
and give a proof of extremality by using the walk counting arguments covered in the previous sections.
In the cases $\ell \in \{2,3,5\}$, the bound (\ref{eq:lv}) for $\mbox{ex}(n,\mc{C}_{2\ell})$ is asymptotically tight,
i.e.\ $\mbox{ex}(n,\mc{C}_{2\ell}) \sim \frac{1}{2}n^{1 + 1/\ell}$.
This is shown by constructions from projective geometry, due to Erd\H{o}s and R\'enyi \cite{ER} for $\ell=2$,
and to Benson \cite{Ben} and Singleton \cite{S} for $\ell=3$ and $\ell=5$, based on generalized polygons.
We briefly describe these here, referring the reader to~\cite{LUW2} for more details.
Suppose $P$ is a set of points, $L$ a set of lines, and $I$ an incidence relation between $P$ and $L$.
The {\em bipartite incidence graph} is the bipartite graph with parts $P$ and $L$
such that $p \in P$ is adjacent to $\ell \in L$ if and only if $(p,\ell) \in I$.
Let $q \geq 2$ and $\ell \geq 2$ be integers.
A {\em generalized $(\ell+1)$-gon of order $q$} consists of a set of $q^{\ell} + q^{\ell - 1} + \dots + 1$ points $P$ and a set of $q^{\ell} + q^{\ell - 1} + \dots + 1$ lines $L$ with an incidence relation $I$ such that the bipartite incidence graph of $(P,L,I)$ is $(q + 1)$-regular and has girth $2\ell+2$ and diameter $\ell+1$. Feit and Higman \cite{FH} showed that
generalized $(\ell+1)$-gons of order $q$ only exist for $\ell = 2,3,5$. These are known as generalized triangles,
generalized quadrangles and generalized hexagons respectively. Generalised triangles are precisely projective planes,
whereas the generalized quadrangles and generalized hexagons were first constructed by Tits in his seminal paper~\cite{T}.

We start with the proof of an upper bound for $\mbox{z}(n,\mc{C}_{2\ell})$, and the characterization of equality by generalized polygons described above. For $\ell=2$ this result is well-known (see \cite{AS}, pp.\ 273--274 for a nice exposition). Here it is convenient to just give the proof for odd $\ell > 2$, as we do not need the other cases.

\begin{prop}\label{even=}
Let $\ell=2$ or $\ell>2$ be odd.
Suppose $n \in \mathbb N$, and let $q \in \mathbb R^+$ be defined by $n = 2(q^{\ell} + q^{\ell-1} +  \cdots + 1)$.
Then $\mbox{z}(n,\mc{C}_{2\ell}) \leq \frac{1}{2}(q + 1)n$, with equality if and only if $q$ is a positive integer and
there exists a generalized $(\ell + 1)$-gon on $n$ points.
\end{prop}

\nib{Proof.}
Let $\ell>2$ be odd and $G$ be a bipartite graph on $n$ vertices with $e(G) \ge (q+1)n/2$
containing no even cycle of length at most $2\ell$.
Then $G$ has average degree at least $q+1$.
Let $A$ and $B$ be the parts of $G$. Then $|A|+|B|=n$.
By the girth assumption, for any ordered pair of vertices $x,y$
there is at most one non-returning walk of length at most $\ell$ from $x$ to $y$.
Since $\ell$ is odd, we need only consider pairs with
$x$ and $y$ in different parts, of which there are $2|A||B| \le n^2/2$.
By Proposition \ref{nwalks}, the number of non-returning walks of odd length at most $\ell$ is at least
\[n(q+1)\sum_{i=0}^{(\ell-1)/2} q^{2i} = n\sum_{j=0}^\ell q^j = n^2/2.\]
We conclude that the number of non-returning walks of length at most $\ell$ is exactly $2|A||B| = n^2/2$,
so $|A| = |B| = n/2$ and $G$ has diameter $\ell + 1$. Since equality only holds in Proposition \ref{nwalks}
for regular graphs, every vertex of $G$ has degree $q+1$.
Therefore $G$ is the bipartite incidence graph of a generalized $(\ell+1)$-gon. \qed

We also require the following bound on unbalanced Zarankiewicz numbers $\mbox{z}(a,b,\mc{C}_{2\ell})$.

\begin{prop} \label{unbalanced}
Suppose $\ell=2$ or $\ell \ge 3$ is odd, and $a,b \geq 1$. Then
\[ \mbox{z}(a,b,\mc{C}_{2\ell}) \le (ab)^{\frac{1}{2} + \frac{1}{2\ell}} + \max\{a,b\}.\]
\end{prop}

\nib{Proof.}
The case $\ell = 2$ follows from a slightly stronger result of Reiman~\cite{R}.
Now suppose $\ell \geq 3$ is odd.
Let $G$ be a $\mc{C}_{2\ell}$-free bipartite graph with parts $A$ and $B$ of sizes $a$ and $b$,
and let $\alpha$ and $\beta$ denote the average degrees of vertices in $A$ and $B$.
As in Proposition \ref{even=}, there are at most $2ab$ non-returning walks of odd length at most $\ell$.
Using the lower bound from Proposition \ref{bipnoreturn}, we deduce that
\begin{equation}\label{eq:unbalanced}
2e(G) \sum_{i=0}^{(\ell - 1)/2} (\alpha - 1)^i(\beta - 1)^i \le 2ab.
\end{equation}
Suppose for a contradiction that $e(G) > (ab)^c + \max\{a,b\}$, where $c=\frac{1}{2} + \frac{1}{2\ell}$.
Then $\alpha-1=(e(G)-a)/a > a^{c-1}b^c$ and $\beta-1=(e(G)-b)/b > a^cb^{c-1}$, so
$(\alpha - 1)(\beta - 1) > (ab)^{2c-1}=(ab)^{1/\ell}$. However, this gives
\[e(G)(\alpha - 1)^{(\ell - 1)/2}(\beta - 1)^{(\ell - 1)/2} > (ab)^c (ab)^{(\ell-1)/2\ell} = ab,\]
which contradicts (\ref{eq:unbalanced}). This gives the required bound on $e(G)$. \qed

\section{Degrees in nearly extremal graphs}

In this section, we show that in nearly extremal $\mc{C}_{2\ell}$-free bipartite graphs for $\ell \in \{2,3,5\}$, the
number of edges containing a vertex of degree substantially more than the average degree is small.
This is used in Step 1 of the stability result (Controlling the maximum degree)
and Step 1 of the exact result (Pseudorandomness).
We also show how to classify extremal graphs once extremal graphs of large minimum degree are classified.
This is used in Step 2 of the exact result (Controlling the minimum degree).

\subsection{Bounding the maximum degree of $C_4$-free graphs}

First we show that any $C_4$-free bipartite graph on $n$ vertices does not contain many edges
on vertices of degree much more than roughly $(n/2)^{1/2}$, which is the average degree in extremal
$C_4$-free bipartite graphs by Proposition \ref{even=}.
Recall that we write $d_B(v)$ for the number of neighbors of a vertex $v$ in a set $B$
and $e(S,B)$ for the number of ordered pairs $(s,b)$ with $s \in S$, $b \in B$ and $sb \in E(G)$.

\begin{lem}\label{c4maxdeg}
Suppose $G$ is a $C_4$-free graph on $n$ vertices and $A,B \subset V(G)$. Suppose $0 < \eps < \sqrt{3}$
and let \[ S = \{v \in A : d_B(v) \geq (1+\eps)|B|^{1/2}\}.\]
Then $e(S,B) \le 2|B|/\eps$. In particular, we have the following:
\begin{itemize}
\item[(i)] If $G$ is a $C_4$-free graph then at most $2n/\eps$ edges contain vertices of degree at least $(1+\eps)\sqrt{n}$.
\item[(ii)] If $G$ is a $C_4$-free bipartite graph with parts $X$ and $Y$, then
$e(T,X \cup Y) \leq 2n/\eps$, where
\[ T = \{v \in X : d_Y(v) \geq (1 + \eps)|Y|^{1/2}\} \cup \{v \in Y : d_X(v) \geq (1 + \eps)|X|^{1/2}\}.\]
\end{itemize}
\end{lem}

\nib{Proof.}
Since $G$ is $C_4$-free, every choice of a vertex $v$ in $B$ and two neighbors $s,s'$ of $v$ in $S$
gives rise to a different pair $\{s,s'\}$ in $S$, so
\[ \sum_{v \in B} \binom{d_S(v)}{2} \le \binom{|S|}{2}.\]
Since the function $f(x)=x(x-1)$ is convex for $x \ge 1$,
Jensen's inequality gives
\[ |S|(|S|-1) \ge e(S,B)(e(S,B)/|B|-1).\]
Suppose for a contradiction that $e(S,B) > 2|B|/\eps$. Then
\[ e(S,B) - |B| >(1-\eps/2)e(S,B).\]
Since $e(S,B) \ge |S| \cdot (1+\eps)|B|^{1/2}$ by definition of $S$, we have
\[ |B||S|^2 \ge e(S,B)(e(S,B)-|B|) \ge (1-\eps/2)e(S,B)^2 \ge (1-\eps/2)(1+\eps)^2|B||S|^2.\]
However, $(1-\eps/2)(1+\eps)^2 = 1 + (3-\eps^2)\eps/2 > 1$, contradiction.
So we must instead have $e(S,B) \le 2|B|/\eps$.
Now statement (i) follows by taking $A=B=V(G)$, and statement (ii) follows by first taking $A=X$ and $B=Y$ and then
taking $A=Y$ and $B=X$.
\qed

\subsection{Bounding the maximum degree in $\mc{C}_{2\ell}$-free graphs}

Next we prove an analogue of Lemma \ref{c4maxdeg} for $\mc{C}_{2\ell}$-free bipartite graphs.
The accurate estimate given here is only needed for the pseudorandomness argument;
a cruder ad hoc method will suffice for bounding the maximum degree in Theorem \ref{stability}.
In the following lemma we could give a quantitative description of
how the $o(\cdot)$ estimate depends on $d$; however the calculations are somewhat heavy,
so we rather retain only the asymptotic statement as $n\to\infty$.

\begin{lem}\label{evenmaxdeg}
Let $\eps > 0$ and suppose $G'$ is a $\mc{C}_{2\ell}$-free bipartite graph
on $n$ vertices with average degree $d \sim (n/2)^{1/\ell}$, where $\ell$ is odd.
Then at most $o(e(G'))$ edges contain a vertex of degree at least $(1 + \eps)d$.
\end{lem}

\nib{Proof.}
We start by repeatedly removing vertices of degree $0$ or $1$ to leave a graph $G$ with minimum degree at least $2$.
This process removes at most $n=o(e(G'))$ edges, so Proposition \ref{even=} implies that we remove $o(n)$ vertices.
Let $A$ and $B$ be the parts of $G$. As in the proof of Propositions \ref{even=} and \ref{unbalanced},
the girth assumption implies that the number of non-returning walks of odd length at most $\ell$
is at most $2|A||B|$. Here we will just use the estimate
\begin{equation}\label{eq:upper}
\nu_\ell(G) \le 2|A||B|/n \le n/2.
\end{equation}
In the lower bound for $\nu_\ell$ we will use the full strength
of Proposition \ref{bipnoreturn} to get an improvement
if there are many edges incident to vertices of large degree.
First we show that the bipartition is roughly balanced.
Proposition \ref{unbalanced} gives $e(G) \leq (|A||B|)^{\frac{1}{2} + \frac{1}{2\ell}} + n$.
Since $e(G) \sim e(G') = nd/2 \sim (n/2)^{1+1/\ell}$ we deduce that $|A| \sim |B| \sim n/2$.

Now let $S$ be the set of vertices in $A$ of degree at least $(1 + \eps)d$.
We will show that $e(S,B) = o(e(G))$. Proposition \ref{bipnoreturn} gives
\begin{equation}\label{eq:lower}
\nu_{\ell}(G) \geq d \pi(A)^{(\ell - 1)/2} \pi(B)^{(\ell - 1)/2},
\end{equation}
where for any $C \sub V(G)$ we write
\[\pi(C) = \prod_{v \in C} (d(v) - 1)^{d(v)/e(G)}.\]
For $\pi(B)$ we just use the simple bound
\[ \pi(B) \ge e(G)/|B| -1 \sim d.\]
To see this we apply Jensen's inequality with the convex function $f(x)=x\log(x-1)$ for $x \ge 2$,
using the fact that $G$ has minimum degree at least $2$. This shows that $\pi(B)$
is minimized when $d(v)= e(G)/|B|$ for all $v \in B$.
Since $e(G) \sim (n/2)^{1+1/\ell}$ and $|B| \sim n/2$ we get the stated bound on $\pi(B)$.

For $\pi(A)$ we estimate $\pi(S)$ and $\pi(A \sm S)$ separately. We write
\[ |S| = \sigma n/2 \ \mbox{ and } e(S,B) = \rho e(G).\]
The parameters $\rho$ and $\sigma$ satisfy
\begin{equation}\label{eq:rho}
(1+o(1))(1+\eps)\sigma \le \rho \le (1+o(1))\sigma^{\frac{1}{2} + \frac{1}{2\ell}},
\end{equation}
where the lower bound follows from $e(S,B) \ge |S|(1+\eps)d$, and the upper bound from
Proposition \ref{unbalanced}, which gives $e(S,B) \le (|S||B|)^{\frac{1}{2} + \frac{1}{2\ell}} + n$.
By Jensen's inequality, $\pi(S)$ is minimized when $d(v)=\rho e(G)/|S|$ for all $v \in S$,
and $\pi(A \sm S)$ is minimized when $d(v)=(1-\rho)e(G)/|A\sm S|$ for all $v \in A \sm S$.
Therefore
\begin{align*}
\pi(A)=\pi(S)\pi(A \sm S) \ge \brac{\frac{\rho e(G)}{|S|}-1}^\rho \brac{\frac{(1-\rho)e(G)}{|A\sm S|}-1}^{1-\rho}
\sim d \brac{\frac{\rho}{\sigma}}^\rho \brac{\frac{1-\rho}{1-\sigma}}^{1-\rho}.
\end{align*}
Applying (\ref{eq:upper}) and (\ref{eq:lower}) we deduce that
\[n/2 \ge \nu_\ell(G) \ge (1+o(1)) d
\brac{d^2\brac{\frac{\rho}{\sigma}}^\rho \brac{\frac{1-\rho}{1-\sigma}}^{1-\rho}}^{(\ell - 1)/2}.\]
Since $d \sim (n/2)^{1/\ell}$ we obtain
\[ \brac{\frac{\rho}{\sigma}}^\rho \brac{\frac{1-\rho}{1-\sigma}}^{1-\rho} \le 1 + o(1). \]
Now recall that our goal is to show that $e(S,B) = o(e(G))$, i.e.\ $\rho=o(1)$.
Suppose for a contradiction that we can choose an infinite sequence of graphs
$G_n$ as above with analogous parameters $\rho_n$ and $\sigma_n$ such that $\rho_n\ge c$
for some constant $c>0$. Then we can pass to a subsequence such that $\sigma_n \to s$
and $\rho_n \to r$, for some $r,s \in [0,1]$ with $r \ge c$.
Then $r$ and $s$ satisfy
\begin{equation}\label{eq:frs}
f(r,s) := \brac{\frac{r}{s}}^r \brac{\frac{1-r}{1-s}}^{1-r} \le 1.
\end{equation}
We also have $(1+\eps)s \le r \le s^{\frac{1}{2} + \frac{1}{2\ell}}$ by (\ref{eq:rho}).
Note that this implies $r<1$.
Consider the function $g(t)=f(r,t)$. Computation of derivatives gives
\[g'(t) = \frac{t-r}{t(1-t)}g(t) \ \mbox{ and }\
g''(t) = \frac{2(t-r)^2+r(1-r)}{t^2(1-t)^2} g(t).\]
Thus  $g'(r)=0$ and $g''(t)>0$ for all $t$, so
$g(t)$ is convex, is decreasing for $t \leq r$  and is minimized at $t=r$.
Since $s \le r/(1+\eps)$ the minimum value possible for $f(r,s)$ is at $s=r/(1+\eps)$.
Substituting in (\ref{eq:frs}) and simplifying gives
\[ 1+\eps \le \brac{1+\frac{\eps}{1-r}}^{1-r}.\]
However this is a contradiction, by the standard inequality $(1+y/x)^x<1+y$ for $x,y \in (0,1)$.
We deduce that $e(S,B) = o(e(G))$.
The same argument shows that $o(e(G))$ edges contain a vertex in $B$
of degree at least $(1+\eps)d$, so the proof is complete.
\qed

\subsection{Bounding the minimum degree}\label{sec:mindeg}

Here we implement Step 2 in the proof of the exact result,
by reducing the proof of Theorem \ref{even} to the case when the minimum degree satisfies $\delta(G)> q/4$.
Recall that $G$ is a $\mc{C}_{2\ell}$-free graph on $n$ vertices with at least $\frac{1}{2}(q + 1)n$ edges, where
$n$ is large, $\ell \in \{2,3,5\}$ and $q \in \mathbb R^+$ is defined by $n=2(q^\ell + q^{\ell - 1} + \dots + 1)$.
We use the vertex deletion argument as in Section \ref{sec:c4}.
Consider a sequence of graphs $G=G_n,G_{n-1},\cdots,G_t$ for some $0 \le t \le n$,
where $G_{i-1}$ is obtained from $G_i$ by deleting a vertex of degree at most $q/4$, while possible.
Note that, by definition, the minimum degree of $G_t$ is at least $q/4$.
We claim that $t>n/2$. For suppose otherwise and consider $G_{n/2}$.
We have $e(G_{n/2}) \le \frac{1}{2}(n/2)^{1 + 1/\ell} + O(n)$ by (\ref{eq:lv}) (see Introduction).
Also, the number of edges deleted is at most $(n/2)(q/4)$. Since we assume
that $e(G) \ge (q+1)n/2$ we get $e(G_{n/2}) > \frac{3}{4} qn/2$.
But $qn/2 \sim (n/2)^{1+1/\ell}$, so this contradicts the upper bound.
Thus we do have $t>n/2$. Next we claim that $e(G_t) \ge (r+1)t/2$,
where $r$ is defined by $t=2\sum_{i=0}^\ell r^i$. To see this note
that $e(G_t) \ge e(G) - (n-t)q/4 \ge (q+1)n/2 - (n-t)q/4$. Then, using $q \ge r$, we calculate
\begin{eqnarray*}
(q+1)n/2-(r+1)t/2-(n-t)q/4
 &=& (q+1)\sum_{i=0}^\ell q^i - (r+1)\sum_{i=0}^\ell r^i - (q/2)\brac{\sum_{i=0}^\ell q^i-\sum_{i=0}^\ell r^i} \\
&=& (q/2+1)\sum_{i=0}^\ell q^i - (r+1-q/2)\sum_{i=0}^\ell r^i \\
&\ge& (q/2+1)\brac{\sum_{i=0}^\ell q^i - \sum_{i=0}^\ell r^i} \; \; \ge \; \; 0.
\end{eqnarray*}
We deduce that $e(G_t) \ge (r+1)t/2$, with strict inequality unless $t=n$.
Now suppose that we know Theorem \ref{even} holds
under the additional assumption $\delta(G)> q/4$.
Applying this to $G_t$ gives $e(G_t) \le (r+1)t/2$.
Thus we must have $t=n$, so $\delta(G)>q/4$,
and we are justified in assuming this when proving Theorem \ref{even}.

\section{Pseudorandomness}

A key ingredient in the proofs of  Theorem \ref{stability} and \ref{even} is the notion of {\em pseudorandomness}.
There are many equivalent notions of pseudorandomness in graphs: we refer the reader to \cite{KS} for a survey.
In this section we will present spectral properties of graphs that imply pseudorandomness.
We will prove the following result, which expresses the pseudorandomness property of a
$\mc{C}_{2\ell}$-free bipartite graph of close to maximum size:
for any two large sets the number of edges between them
is roughly the same as in a random bipartite graph of the same density.

\begin{thm}\label{bippr}
Suppose $G$ is a $\mc{C}_{2\ell}$-free bipartite graph on $n$ vertices with
parts $X$ and $Y$ and average degree $d \sim (n/2)^{1/\ell}$.
Then for any $S \sub X$ and $T \sub Y$ we have $e(S,T) = \frac{2d}{n}|S||T| + o(n^{1+1/\ell})$.
\end{thm}

This theorem for $\ell = 2$ is the pseudorandomness part of the sketch proof given in Section \ref{sketch-exact}.

\subsection{Pseudorandomness of regular graphs}

The exposition in this subsection repeats that in \cite[Section 2.4]{KS},
so we will be brief and refer the reader to that survey for more details.
As a warmup, we give an exposition of the fact, first proved by N. Alon, that a regular graph with a large spectral
gap is pseudorandom (this is sometimes known as the `expander mixing lemma').
Then we establish analogous results in the bipartite setting,
which seem not have been explicitly presented in the previous literature.

Suppose $G$ is a graph on $n$ vertices and let $A =(a_{uv})_{u,v \in V(G)}$
be its adjacency matrix, i.e.\ $a_{uv}$ is $1$ if $uv$ is an edge or $0$ otherwise.
(We fix some ordering $v_1,v_2,\dots,v_n$ of the vertices and identify $v_i$ with $i$.)
Then $A$ is a real symmetric matrix, so has an orthonormal basis $x_1,x_2,\dots,x_n$
of eigenvectors with real eigenvalues, which we order so that $\lambda_1 \ge \lambda_2 \ge \dots \ge \lambda_n$.
Note that $\sum_{i=1}^n \lambda_i = \mbox{Tr}(A)=0$, so $\lambda_1 \ge 0$ and $\lambda_n \le 0$.
The Perron-Frobenius theorem implies that $|\lambda_i| \le \lambda_1$ for all $i$
and all entries of $x_1$ are non-negative.
If $G$ is $d$-regular then we have $x_1 = e := n^{-1/2}(1,1,\dots,1)$ and $\lambda_1=d$.
In this case it is easy to verify that $|\lambda_i| \le d$ for all $i$,
since if $x_{i,j}$ has the largest absolute value among the coordinates of $x_i$ then
\[ |\lambda_i x_{i,j}|=|(Ax_i)_j| = \bsize{\sum_{k \in N(j)} x_{i,k}} \le d|x_{i,j}|.\]
In this non-bipartite setting we write
\[ \lambda = \max_{i \ne 1} |\lambda_i|.\]
We have the following pseudorandomness property for regular graphs (see, e.g., \cite{AS, KS}), whose short proof we include
for the convenience of the reader.

\begin{lem}\label{pseudo-reg-graph}
Suppose $G$ is a $d$-regular graph on $n$ vertices.
Then for any $S,T \sub V(G)$ we have
\[ \bsize{e(S,T) - \frac{d}{n}|S||T|} \leq \lambda\sqrt{|S||T|}.\]
\end{lem}

\nib{Proof.}
Let $\chi_S$ and $\chi_T$ denote the characteristic vectors of $S$ and $T$,
which are equal to $1$ or $0$ in position $v$ according as
$v$ belongs or does not belong to the corresponding set.
Then $e(S,T) = \chi_S^t A \chi_T$. Let $\{x_1,x_2,\dots,x_n\}$ be an orthonormal basis of eigenvectors,
where $x_i$ is the eigenvector corresponding to eigenvalue $\lambda_i$, and write
$\chi_S = \sum_{i=1}^n s_i x_i$ and $\chi_T = \sum_{i=1}^n t_i x_i$. Then $e(S,T) = \sum_{i=1}^n \lambda_i s_i t_i$.
Note that $\sum_{i=1}^n s_i^2 = \chi_S \cdot \chi_S = |S|$
and similarly $\sum_{i=1}^n t_i^2 = |T|$. Thus we can estimate
\[ \bsize{\sum_{i > 1} \lambda_i s_i t_i} \le \lambda \sum_{i > 1} |s_i||t_i|\le \lambda \sqrt{|S||T|}\]
by the Cauchy-Schwarz inequality.
Since $\lambda_1=d$, $s_1 = e \cdot \chi_S = n^{-1/2}|S|$, and similarly $t_1 = n^{-1/2}|T|$,
we obtain
\begin{equation}\tag*{$\Box$}
\bsize{e(S,T) - \frac{d}{n}|S||T|} \leq \lambda\sqrt{|S||T|}.
\end{equation}

\subsection{Pseudorandomness of regular bipartite graphs}

Now we adapt the arguments of the previous subsection to the bipartite setting.
Let $G$ be a bipartite graph on $n$ vertices with parts $X$ and $Y$.
We choose the ordering $v_1,v_2,\dots,v_n$ of the vertices so that $X$ precedes $Y$.
When we consider vectors of length $n$ this gives us a natural partition of its
co-ordinates into two blocks corresponding to $X$ and $Y$.
Then the adjacency matrix $A$ has block structure $\binom{0\ M}{M^t\ 0}$
where $M$ is the bipartite incidence matrix of $G$,
i.e.\ $M$ has rows indexed by $X$, columns indexed by $Y$,
and entries $m_{xy}$ equal to $1$ if $xy$ is an edge, otherwise $0$.
First we consider the case when $G$ is $d$-regular, which implies that $|X|=|Y|=n/2$.
Then we have $\lambda_1=d$, with eigenvector $e:=n^{-1/2}(1,\dots,1)$ as before,
and $\lambda_n=-d$, with eigenvector $\ov{e}:=n^{-1/2}(1,\dots,1,-1,\dots,-1)$
having $n^{-1/2}$ in its $X$-coordinates and $-n^{-1/2}$ in its $Y$-coordinates.
In the bipartite setting we re-define $\lambda$ by
\[ \lambda = \max_{i \ne 1,n} |\lambda_i|.\]
We have the following pseudorandomness property for regular bipartite graphs.

\begin{lem} \label{pseudo-reg-bip}
Suppose $G$ is a $d$-regular bipartite graph on $n$ vertices with parts $X$ and $Y$. Then for any $S \sub X$ and $T \sub Y$
we have
\[ \bsize{e(S,T) - \frac{2d}{n}|S||T|} \leq \lambda\frac{|S|+|T|}{2}.\]
\end{lem}

\nib{Proof.}
Let $\chi=(\chi_S,\chi_T)$ denote the characteristic vector of $S \cup T$.
Then
\[ \chi^t A \chi = \chi_S^t M \chi_T + \chi_T^t M^t \chi_S = 2e(S,T).\]
Writing $\chi = \sum_{i=1}^n a_i x_i$ in the eigenvector basis
we obtain
\[ 2e(S,T) = \sum_{i=1}^n \lambda_i a_i^2.\]
Since $\sum_{i=1}^n a_i^2 = \chi \cdot \chi = |S|+|T|$
we can estimate
\[ \bsize{\sum_{i \ne 1,n} \lambda_i a_i^2} \le |\lambda|(|S|+|T|).\]
We have
\[ \lambda_1=d,\quad a_1 = e \cdot \chi = n^{-1/2}(|S|+|T|), \quad
\lambda_n=-d \quad \mbox{and} \quad a_n=\ov{e} \cdot \chi = n^{-1/2}(|S|-|T|),\]
so $\lambda_1 a_1^2 + \lambda_n a_n^2 = \frac{4d}{n}|S||T|$.
This gives the stated estimate for $e(S,T)$.
\qed

\bigskip

\nib{Remark.} An alternative derivation of similar estimates may be obtained
from the singular value decomposition of the bipartite incidence matrix.

\subsection{Nearly regular bipartite graphs}

We want to show pseudorandomness for $\mc{C}_{2\ell} \cup \{C_k\}$-free graphs which are nearly extremal.
Such graphs are not necessarily regular -- for instance they may contain isolated vertices -- so to treat them we
will prove an analogue of Lemma \ref{pseudo-reg-bip} for nearly regular bipartite graphs.
The quantity we use to measure irregularity of a graph is its {\em variance}. If $G$ is a graph of average degree $d$
with $n$ vertices, then the variance of $G$ is defined by
\[ \mbox{\sc Var}(G) : = \frac{1}{n} \sum_v (d(v)-d)^2 = \frac{1}{n} \sum_v d(v)^2-d^2.\]
As before we set $\lambda = \lambda(G) = \max_{i \ne 1,n} |\lambda_i|$ in the bipartite setting.
We have the following pseudorandomness property for nearly regular bipartite graphs
(for a similar statement for non-bipartite graphs see \cite{KS}).

\begin{lem}  \label{pseudo-near-reg-bip}
Let $\beta,\gamma \in (0,1)$ and $\alpha = 4\beta^{1/2}\gamma^{-1}<1/4$.
Suppose $G = G(X,Y)$ is a bipartite graph on $n$ vertices with average degree $d$ and
\begin{itemize}
\item[(i)]  $\lambda(G) < (1 - \gamma)d$,
\item[(ii)] $\mbox{\sc Var}(G) < \beta d^2$.
\end{itemize}
Then for any $S \sub X$ and $T \sub Y$ we have
\[ \bsize{e(S,T) - \frac{2d}{n}|S||T|} \leq (4\alpha d + \lambda/2)n.\]
\end{lem}

Recall that $e=n^{-1/2}(1,\dots,1)$, and $\ov{e}:=n^{-1/2}(1,\dots,1,-1,\dots,-1)$
has $n^{-1/2}$ in its $X$-coordinates and $-n^{-1/2}$ in its $Y$-coordinates.
The following estimates will be used in the proof and later in the paper.
\begin{align*}
\lambda_1 & = \max\{ x^t A x: \|x\|=1 \} \ge e^t A e = n^{-1}\sum_v d(v) = d, \quad \mbox{ and }\\
\lambda_n & = \min\{ x^t A x: \|x\|=1 \} \le \ov{e}^t A \ov{e} = n^{-1}\sum_v -d(v) = -d.
\end{align*}

\nib{Proof.}
Write $|S|=s$ and $|T|=t$. Without loss of generality $s \ge t$. As in Lemma \ref{pseudo-reg-bip},
we consider the characteristic vector $\chi=(\chi_S,\chi_T)$ and write $\chi = \sum_{i=1}^n a_i x_i$, where $\{x_i : 1 \leq i \leq n\}$ is an orthonormal basis
of eigenvectors and $x_i$ is the eigenvector for $\lambda_i$.
Then $\sum_{i=1}^n a_i^2=\chi^t \cdot \chi=s+t$  $2e(S,T)= \chi^t A \chi = \sum_{i=1}^n \lambda_i a_i^2$,  and we estimate, with $\lambda = \lambda(G)$:
\begin{equation}\label{eq:rest}
|2e(S,T) - \lambda_1 a_1^2 - \lambda_n a_n^2|
=\bsize{\sum_{i \ne 1,n} \lambda_i a_i^2} \le \lambda(s+t) \le \lambda n.
\end{equation}
Next write $e = \sum_{i=1}^n c_i x_i$ and $\ov{e} = \sum_{i=1}^n \ov{c}_i x_i$,
where $\sum c_i^2 = \sum \ov{c}_i^2 = 1$; since we can replace any eigenvector $x_i$
by $-x_i$ we can assume that $c_1>0$ and $\ov{c}_n>0$.
We will show that $x_1$ is close to $e$ and $x_n$ is close to $\ov{e}$.
Consider $z=Ae-de$, which has co-ordinates $z_i = n^{-1/2}(d(v_i)-d)$,
and $\ov{z}=A\ov{e}+d\ov{e}$, which has co-ordinates $\ov{z}_i= \pm n^{-1/2}(d(v_i)-d)$,
with positive sign when $v_i \in Y$ and negative sign when $v_i \in X$.
Then by (ii), \[ z \cdot z = \ov{z} \cdot \ov{z} = n^{-1} \sum_v (d(v)-d)^2 = \mbox{\sc Var}(G) \leq \beta d^2.\]
We also have $z=\sum (\lambda_i-d)c_ix_i$ and $\ov{z}=\sum (\lambda_i+d)\ov{c}_ix_i$,
so \[ z \cdot z = \ov{z} \cdot \ov{z} = \sum c_i^2(\lambda_i-d)^2 = \sum \ov{c}_i^2(\lambda_i+d)^2.\]
Then $\sum_{i \ne 1} c_i^2 \le (d-\lambda)^{-2} \beta d^2 \le \beta \gamma^{-2}$ by (i),
so $c_1 \ge c_1^2 \ge 1 - \beta \gamma^{-2}$,
and
$$\|e-x_1\|^2 = (1-c_1)^2 + \sum_{i \ne 1} c_i^2 = (1-c_1)^2+1-c_1^2=2(1-c_1) \le 2\beta \gamma^{-2}.$$
Similarly, one can prove that  $\ov{c}_n \geq 1 - \beta \gamma^{-2}$ and $\|\ov{e}-x_n\|^2 \le 2\beta \gamma^{-2}$. Next we estimate $\lambda_1$ and $\lambda_n$.
Recall that $\lambda_1 = \|Ax_1\|$ and $d = \|de\|$, and so by the triangle inequality,
\begin{equation}
\label{eq5.4}
 |\lambda_1 - d| = \big|\|Ax_1\| - \|de\|\big| \leq \| Ax_1 - de \| \leq \|Ae-de\| + \|Ae - Ax_1\|.
\end{equation}
We noted before the proof that $\lambda_1 \ge d$, so we can write $\lambda_1 = (1+q_1)d$ with $q_1 \ge 0$.
Since $\|Ae-de\|=\|z\| < \beta^{1/2}d$ and $\|A(e-x_1)\| \le \lambda_1 \|e-x_1\| \le (2\beta)^{1/2}\gamma^{-1} \lambda_1$,
from (\ref{eq5.4}) we have $q_1 \le \beta^{1/2} + (2\beta)^{1/2}\gamma^{-1} (1+q_1)$. Thus, recalling that $\alpha=4\beta^{1/2}\gamma^{-1}<1/4$, we can estimate
$$q_1 \leq \frac{\beta^{1/2} + (2\beta)^{1/2}\gamma^{-1}}{1-(2\beta)^{1/2}\gamma^{-1}}\leq \frac{\alpha/4+\alpha/(2\sqrt{2})}{1-1/(8\sqrt{2})}
< \alpha,$$
i.e.\ \[  d \le \lambda_1 < (1 + \alpha)d.\]
Similarly, we have \[  -d \ge \lambda_n > -(1 + \alpha)d.\]
We also estimate $a_1 = \chi \cdot x_1$ by $\chi \cdot e = n^{-1/2}(s+t)$ and the inequality
\[ \|\chi \cdot (e-x_1)\|^2 \le \|\chi\|^2\|e-x_1\|^2 \le 2\beta \gamma^{-2}(s+t) < \alpha^2 (s + t),\]
which gives
\[ |a_1 - n^{-1/2}(s+t)| \leq \alpha (s+t)^{1/2}.\]
Similarly, we estimate $a_n = \chi \cdot x_n$ by $\chi \cdot \ov{e} = n^{-1/2}(s-t)$, and get
\[ |a_n - n^{-1/2}(s-t)| \leq \alpha (s+t)^{1/2}.\]
Now we have the necessary ingredients to estimate $\lambda_1 a_1^2 + \lambda_n a_n^2$.
Recall that $\lambda_1 \geq d$, $\lambda_n \geq -d(1 + \alpha)$ and $s \ge t$.
Using the above estimates for $a_1, a_n$, we have the lower bounds
\begin{align*}
\lambda_1 a_1^2 & \ge d (n^{-1/2}(s+t) - \alpha (s+t)^{1/2})^2, \mbox{ and } \\
\lambda_n a_n^2 & \ge -\; d(1 + \alpha) (n^{-1/2}(s-t) + \alpha  (s+t)^{1/2})^2.
\end{align*}
Thus we obtain
\[ \lambda_1 a_1^2 + \lambda_n a_n^2 \ge \frac{4d}{n}st - \frac{\alpha d}{n} (s-t)^2
- 2\alpha d\brac{\frac{s+t}{n}}^{1/2}\Big((s+t)+(1+\alpha)(s-t)\Big) - \alpha^3 d(s+t).\]
Since $s+t \le n$ and $\alpha<1/2$ we get
\[ \lambda_1 a_1^2 + \lambda_n a_n^2 - \frac{4d}{n}st \ge - 8\alpha dn.\]
The estimates for the upper bound are similar but slightly more technical. We use
\begin{align*}
\lambda_1 a_1^2 & \leq d(1 + \alpha)(n^{-1/2}(s+t) + \alpha (s+t)^{1/2})^2, \mbox{ and } \\
\lambda_n a_n^2 & \leq
\begin{cases}
-\; d(n^{-1/2}(s-t) - \alpha  (s+t)^{1/2})^2 & \mbox{ if } \alpha  (s+t)^{1/2} \le n^{-1/2}(s-t), \\
0 & \mbox{ if } \alpha  (s+t)^{1/2} > n^{-1/2}(s-t).
\end{cases}
\end{align*}
In the case $\alpha  (s+t)^{1/2} \le n^{-1/2}(s-t)$ we have
\[ \lambda_1 a_1^2 + \lambda_n a_n^2 \le \frac{4d}{n}st + \frac{\alpha d}{n} (s+t)^2
+ 2\alpha d\brac{\frac{s+t}{n}}^{1/2}\Big((1+\alpha)(s+t)+(s-t)\Big) + \alpha^3 d(s+t).\]
In the case $\alpha  (s+t)^{1/2} \ge n^{-1/2}(s-t)$ we have
$(s+t)^2 - 4st = (s-t)^2 \le \alpha^2 n (s+t) \le \alpha^2 n^2$, so
\[ \lambda_1 a_1^2 + \lambda_n a_n^2 \le \lambda_1 a_1^2 \le
d(1 + \alpha)(n^{-1}(4st+\alpha^2 n^2)+ 2\alpha n^{-1/2}(s+t)^{3/2}+\alpha^2(s+t)).\]
Since $s+t \leq n$ and $ \alpha < 1/4$, in both cases we obtain \[ \lambda_1 a_1^2 + \lambda_n a_n^2 - \frac{4d}{n}st \le 8\alpha dn.\]
Combining this with (\ref{eq:rest}) we obtain the stated estimate for $e(S,T)$.
\qed

\subsection{Proof of Theorem \ref{bippr}}

The idea of the proof is to use the connection between eigenvalues and closed walks.
We can control the maximum degree by deleting few edges; then the main contribution
to the upper bound on closed walks of length $2\ell+2$ in Lemma \ref{closedwalks} is $(1/2+o(1))n^2\Delta^2$.
This is very close to $\lambda_1^{2\ell + 2} + \lambda_n^{2\ell + 2} \sim 2d^{2\ell+2}$,
so the other eigenvalues of this graph are small.

We now give the details. Suppose $\eps>0$.
Consider any $\mc{C}_{2\ell}$-free bipartite graph $H$ on $n$ vertices
with parts $X$ and $Y$ and average degree $d \sim (n/2)^{1/\ell}$.
Suppose $S \sub X$ and $T \sub Y$.
By Lemma \ref{evenmaxdeg} there are $o(e(H))$ edges incident to vertices of degree at least $(1 + \eps)d$.
We remove these edges to obtain a graph $G$ of maximum degree $\Delta \leq (1 + \eps)d$
and average degree $d \sim (n/2)^{1/\ell}$. Also, as in the proof of Lemma \ref{evenmaxdeg},
Proposition \ref{unbalanced} gives $|X| \sim |Y| \sim \frac{n}{2}$.
By Proposition \ref{trace} and Lemma \ref{closedwalks}
\[ \frac{1}{n}\sum_{i=1}^n \lambda_i^{2\ell + 2} = w^{\circ}_{2\ell + 2}(G) < (1/2+o(1))n\Delta^2 + (4\Delta)^{\ell + 1}.\]
Now we use the estimates $\lambda_1 \geq d$ and $\lambda_n \leq -d$.
Recalling that $\lambda = \max_{i \ne 1,n} |\lambda_i|$ we have
\[\lambda^{2\ell + 2} \le nw^{\circ}_{2\ell + 2}(G) - 2d^{2\ell + 2}
< (1/2+o(1))n^2\Delta^2 - 2d^{2\ell + 2} + n(4\Delta)^{\ell + 1}.\]
Since $d \sim (n/2)^{1/\ell}$ and $\Delta \leq (1 + \eps)d$
this gives $\lambda^{2\ell + 2} \le ((1+\eps)^2-1+o(1))2d^{2\ell+2}$.
It follows that
\[ \lambda \leq (6\eps)^{1/(2\ell + 2)}d + o(d).\]
Also, the variance is bounded as
$$\mbox{\sc Var}(G)  = \frac{\sum_v (d(v)-d)^2}{n}=
\frac{\sum_v d^2(v)}{n} -d \leq \Delta\cdot \frac{\sum_v d(v)}{n}-d
 \le (1+\eps)d^2 - d^2 = \eps d^2.$$
We now apply Lemma \ref{pseudo-near-reg-bip} with $\beta = \eps$
and $\gamma = 1-(6\eps)^{1/(2\ell + 2)} + o(1)$. Recall
that we need $\alpha=4 \beta^{1/2} \gamma^{-1} <1/4$ , which holds if $\eps$ is small.
Thus
\[ \bsize{e_G(S,T) - \frac{2d}{n}|S||T|} \leq (4\alpha d + \lambda/2)n,\]
where $\alpha = 4\beta^{1/2}\gamma^{-1} < 5\eps^{1/2}$ for small $\eps$.
Recalling that $G$ was obtained from $H$ by deleting $o(n^{1+1/\ell})$ edges we have
\[ \bsize{e_H(S,T) - \frac{2d}{n}|S||T|} < \Big(20\eps^{1/2} + (6\eps)^{1/(2\ell + 2)} + o(1)\Big) n^{1+1/\ell}.\]
Since $\eps$ is arbitrary this proves Theorem \ref{bippr}. \qed

\section{Proof of Theorem \ref{stability}}

Suppose that $G$ is a $\mc{C}_{2\ell} \cup \{C_k\}$-free graph with $n$ vertices with average degree $d=\Theta(n^{1/\ell})$,
where $k>2\ell$ is odd. For the first part of Theorem \ref{stability}, we are required to find a bipartite
graph $H \subset G$ such that $e(H) \geq d^{\ell+1}-o(n^{1+1/\ell})$.
Similarly to the sketch given for $4$-cycles, the idea is that we can take $H$ to be the bipartite subgraph
spanned by $N_\ell(v)$ and $N_{\ell+1}(v)$ for a suitable vertex $v$. The first step is to pass to a subgraph with low
maximum degree.

Let $\Delta = n^{1/\ell + c}$, where $c := 1/2\ell^2$.
Let $S$ be the set of vertices of degree more than $\Delta$,
and let $G_0$ be the graph obtained by removing all edges of $G$ containing at least one vertex of $S$,
where $G_0$ has average degree $d_0$. We will show that $d_0 \sim d$.
To estimate the number of edges removed, recall from (\ref{eq:lv}) that
$\mbox{ex}(n,\mc{C}_{2\ell}) \leq \frac{1}{2}n^{1 + 1/\ell} + O(n)$,
so
\[ |S| \le \frac{n^{1 + 1/\ell} + O(n)}{\Delta} < m := 2n^{1-c}.\]
It follows that
\[ e(S) \leq \frac{1}{2}m^{1 + 1/\ell} + O(m) < n^{1 + 1/\ell - c}.\]
Also Proposition \ref{unbalanced} gives
\[ e(S,V(G)\sm S) < (mn)^{1/2 + 1/2\ell} + n < n^{1 + 1/\ell - c/2}.\]
In particular, we have $e(G_0) > e(G) - o(n^{1 + 1/\ell})$, and therefore
$d_0 \sim d$. For the remainder of the proof we work in the graph $G_0$,
which has maximum degree at most $\Delta$.

Next, by Lemma \ref{paths}, we can choose a vertex $v$ that is the start of at least
$d_0^{\ell+1} - \ell^2\Delta^\ell = d^{\ell+1} + o(n^{1+1/\ell})$ paths of length $\ell+1$ in $G_0$.
We claim that all but $o(n^{1+1/\ell})$ of these paths reach $N_{\ell+1}(v)$.
Consider a breadth-first search tree $T$ rooted at $v$.
Consider any path $P$ of length $\ell + 1$ that does not reach $N_{\ell+1}(v)$.
Then there is a smallest $i$ such that the $(i+1)$st edge of $P$ does
not go from $N_i(v)$ to $N_{i+1}(v)$. By construction of $T$
this edge must either go from $N_i(v)$ to $N_{i-1}(v)$ or lie within $N_i(v)$.
The first case is impossible, as any edge of $E(G_0) \sm E(T)$ between
$N_{i-1}(v)$ and $N_i(v)$ would create an even cycle of length at most $2\ell$.
For the second case we recall that Lemma \ref{find-odd} implies that $G[N_i(v)]$
has average degree at most $2k$ for any $i \le \ell$, since $k>2\ell$.
Also note that from maximum degree assumption we have that $|N_i(v)| \leq \Delta^i$.
This gives at most $k|N_i(v)| \le k\Delta^i$ choices for the $(i+1)$st edge of $P$.
Let $w$ be the first vertex of $P$ in $N_i(v)$.
The subpath of $P$ from $v$ from $w$ is uniquely determined
(otherwise we would have an even cycle of length at most $2\ell$).
Then we have at most $\Delta$ choices for each of the $\ell - i$ subsequent edges of $P$.
In total, the number of choices for $P$ is at most $k\Delta^\ell = o(n^{1+1/\ell})$, as required.

Each edge between $N_{\ell}(v)$ and $N_{\ell+1}(v)$ is contained in at most
one path of length $\ell + 1$ from $v$ to $N_{\ell+1}(v)$, otherwise we would have an even cycle of length at most $2\ell$.
Thus, taking $H$ to be the bipartite graph of edges between $N_{\ell}(v)$ and $N_{\ell + 1}(v)$,
we have $e(H) = d^{\ell+1} + o(n^{1+1/\ell})$,  as required.

To prove the second part of Theorem \ref{stability}, suppose that $d \geq (1 + o(1))(n/2)^{1/\ell}$. The number of edges in the bipartite graph $H$ constructed above is at least $(1+o(1))nd^{\ell+1}/2 \geq (1+o(1))(n/2)^{1 + 1/\ell}$.
On the other hand, from Proposition \ref{even=} we have $e(H) \le (1+o(1))(n/2)^{1+1/\ell}$. Therefore
$e(H) \sim (n/2)^{1 + 1/\ell}$ and $d \le (1 + o(1))(n/2)^{1/\ell}$. Since $d \ge (1+o(1))(n/2)^{1/\ell}$, we also have
$d \sim (n/2)^{1/\ell}$. So $e(H) \sim e(G)$ and by Proposition \ref{even=}, this shows $e(G) \sim \mbox{z}(n,\mc{C}_{2\ell})$.
This completes the proof. \qed

\section{Proof of Theorem \ref{even}}

Suppose that $n$ is large, and $G$ is a graph on $n$ vertices with $e(G) \ge (q+1)n/2$,
where $q \in \mathbb R^+$ is defined by $n=2(q^\ell + q^{\ell - 1} + \dots + 1)$.
Suppose also that $G$ does not contain an even cycle of length at most $2\ell$,
or $C_k$ for some odd $k > 2\ell+2$. We will show that $e(G) = (q + 1)n/2$ and $G$ is bipartite.
Then Proposition \ref{even=} characterizes equality, namely, $G$ must be the incidence graph of a generalized
polygon. We start by considering the case $\ell \ge 2$ and $k > 2\ell + 2$.

{\bf Case 1 : $\ell \ge 2$ and $k > 2\ell + 2$.} Let $H$ be a bipartite subgraph of $G$ with maximum size. We show $G = H$. By Theorem \ref{stability}, $e(H) \sim e(G) \sim (n/2)^{1+1/\ell}$. Furthermore, maximality of $H$ implies that $\delta(H)> \delta(G)/2$,
as if there were a vertex of degree less than $\delta(G)/2$ in $H$
we could move it to the other part and increase the number of edges in $H$.
By Section 4.3, we can assume $\delta(G) > q/4$, so $\delta(H) > q/8$.
Label the parts of $H$ as $X_0$ and $X_1$, with $X_0 \cup X_1 = V(G)$.
Suppose for a contradiction that $G[X_0]$ contains an edge $\{x,y\}$.
Let $z$ be a neighbor of $y$ in $X_1$.
We greedily construct a sequence of mutually disjoint sets $\{y\}$ and $S_x^i$ and $S_z^i$
for $0 \leq i \leq \ell$, where $S_x^0=\{x\}, S_z^0=\{z\}$,
\[ S_x^i \sub N(S_x^{i-1}) \cap X_{i \text{ mod } 2} \quad \mbox{ and }
\quad S_z^i \sub N(S_z^{i-1}) \cap X_{i+1 \text{ mod } 2}.\]
Note that by definition we have $S_x^i \sub N_i(x)$ and $S_z^i \sub N_i(z)$.
By consideration of breadth first search trees as in the proof of Theorem \ref{stability},
every vertex in $N_i(x)$ has exactly one neighbor in $N_{i-1}(x)$
and $e(N_i(x)) < k|N_i(x)|$ for every $1 \le i \le \ell-1$. Moreover, two distinct vertices of
$N_i(x)$ can not have a common neighbor in $N_{i+1}(x)$. Similar statements hold for $z$.
By the minimum degree assumption in $H$, we have $|N_i(x)|, |N_i(z)| \geq (q/8)^i$ for all $ 0 \leq i \leq \ell$.
Since $q \sim (n/2)^{\ell}$, this allows us greedily to choose disjoint sets $S_x^i, S_z^i$ so that
\[ |S_x^i| \sim |S_z^i| \sim (cn)^{i/\ell},\]
for all $1 \leq i \leq \ell$ and some constant $c>0$.
Now we apply Theorem \ref{bippr} to $H$ with $S = S_x^{\ell}$ and $T = S_z^{\ell}$. This gives
\[ e(S,T) \geq \frac{(1+o(1))2(n/2)^{1/\ell}}{n}|S||T| - o(n^{1 + 1/\ell}) > 2kn.\]
In particular, we can find a path of length $k-2\ell-2$ using only edges between $S$ and $T$.
By construction this can be completed to a cycle of length $k$ in $G$, which is a contradiction.
We deduce that $X_0$ is an independent set in $G$. Similarly $X_1$ is independent, so $G=H$.
This completes the proof for $k > 2\ell + 2$.

For the remainder of the proof we consider the special case $\ell = 2$ and $k = 5$.

{\bf Case 2 : $\ell = 2$ and $k = 5$.} We use a similar vertex deletion argument to that in Section 4.3
to find a subgraph $G_t$ of $G$ with $t \geq \lfloor 0.01n\rfloor$ vertices and $\delta(G_t) > 0.51t^{1/2}$.
Starting with $G = G_n$, we produce a graph $G_i$ with $i$ vertices by deleting
a vertex of $G_{i+1}$ of degree at most $0.51 (i + 1)^{1/2}$ for each $i < n$. After
$t$ steps, the total number of edges deleted is less than
\[ \sum_{i=t+1}^n 0.51i^{1/2} < 0.51 \int_{t+1}^{n+1} x^{1/2}\ dx = \frac{1.02}{3} ((n+1)^{3/2}-(t+1)^{3/2}).\]
Suppose that we fail to find a subgraph of minimum degree more than $0.51 t^{1/2}$ in $n - t = \lceil 0.99n \rceil$ steps.
Then we have a graph $G_t$ with $\lfloor 0.01n \rfloor$ vertices and
\begin{eqnarray*}
e(G_t) &>& (q + 1)n/2 - \frac{1.02}{3} ((n+1)^{3/2}-(t+1)^{3/2})
> \frac{1}{2\sqrt{2}}n^{3/2} - \frac{1.02}{3} ((n+1)^{3/2}-(t+1)^{3/2}) \\
&\geq& \Bigl(\frac{1}{2\sqrt{2}} - \frac{1.02}{3} - o(1)\Bigr)n^{3/2} > 0.01n^{3/2} \geq 10 t^{3/2}.
\end{eqnarray*}
This contradicts Theorem \ref{stability}, provided $t = \lfloor 0.01 n \rfloor$ is large enough. So $G$ has a subgraph $G_t$ with $t$ vertices
and $\delta(G_t) > 0.51 t^{1/2}$ and where $t \geq \lfloor 0.01 n\rfloor$.
For $i \le n$ let $q_i$ be the unique positive real defined by $i = 2(q_i^2 + q_i + 1)$,
so that $q_n = q$. Then $q_i = (\sqrt{2i - 3} - 1)/2$, and if $t < n$ then
\[ e(G_{n - 1}) \geq e(G) - 0.51 n^{1/2} \geq (q + 1)n/2 - 0.51 n^{1/2} > (q_{n - 1} + 1)(n - 1)/2\]
for large enough $n$, using
\begin{align*}
& (q + 1)n/2 - (q_{n - 1} + 1)(n - 1)/2  = (\sqrt{2n - 3} - 1)n/4 - (\sqrt{2n - 5} - 1)(n-1)/4 \\
& = \sqrt{2n - 5}/4 + (\sqrt{2n - 3}-\sqrt{2n - 5})n/4 - 1/4
\sim \frac{3}{4\sqrt{2}} n^{1/2} > 0.51 n^{1/2}.
\end{align*}
Repeating this calculation, we get that
\begin{equation}\label{tbound}
\text{if } t < n \text{ then } e(G_t) > (q_t + 1)t/2,
\end{equation}
provided $t \geq \lfloor 0.01 n\rfloor$ is large enough. We will show that $G_t$ is bipartite.

\bigskip

First we pass to a maximum bipartite subgraph $H$ of $G_t$ as in Case 1, with parts $X_0$ and $X_1$.
By Theorem \ref{stability}, $e(H) \sim e(G_t) \sim (t/2)^{3/2}$ . We claim that no vertex $x \in X_0$
has more than $0.09 t^{1/2}$ neighbors in $X_0$, and similarly for $X_1$. To see this, note that such an $x \in X_0$
also has more than $0.09t^{1/2}$ neighbors in $X_1$ by maximality of $H$.
Then $N_2(x)$ contains $\Theta(n)$ vertices of $X_0$ and $\Theta(n)$ vertices of $X_1$. Choose a set $S$
of $\Theta(n)$ vertices of $N_2(x) \cap X_0$ and a set $T$ of $\Theta(n)$ vertices of $N_2(x) \cap X_1$
such that for each $w \in S$ and $z \in T$, there exist paths of length two from $x$ to $w$
and from $x$ to $z$ which share only the vertex $x$. Then in $H$ we apply Theorem \ref{bippr} (pseudorandomness) to
conclude
\[ e(S,T) \geq \frac{(1+o(1))2(t/2)^{1/2}}{t}|S||T| - o(t^{3/2}).\]
In particular, $e(S,T) \neq 0$ and evidently there is a cycle of length five through $x$ and any edge between $S$ and $T$, a contradiction.
Therefore no vertex has more than $0.09t^{1/2}$ vertices in its own part. It follows that every vertex has degree at least $0.501 t^{1/2}$ in $H$.

We next claim that $|X_0| \sim |X_1| \sim t/2$. First note from Proposition \ref{unbalanced} that
\[ e(H) \leq \mbox{z}(t,C_4) \leq (|X_0||X_1|)^{3/4} + \max\{|X_0|,|X_1|\}.\]
On the other hand, $e(H) \sim e(G) \sim (t/2)^{3/2}$, and so we see
\[ \big(|X_0||X_1|\big)^{3/4} \geq (1 + o(1))(t/2)^{3/2}.\]
Since $|X_0| + |X_1| = t$, and we just observed $|X_0||X_1| \geq (1 + o(1))(t/2)^2$, we conclude $|X_0| \sim |X_1| \sim t/2$.

Now we show $G_t[X_0]$ and $G_t[X_1]$ have no edges, so that $G_t$ is bipartite. Suppose for a contradiction that $G_t$ has an edge $\{x,y\}$ with $x,y \in X_0$. Note that $x$ and $y$ have at most one common neighbour, since $G$ is $C_4$-free.
Let $z$ be this common neighbour if it exists, or an arbitrary vertex otherwise.
Let $S$ be the set of ends of paths of length two in $H$ that start at $x$ and avoid $\{y,z\}$.
Let $T$ be the set of ends of paths of length two in $H$ that start at $y$ and avoid $\{x,z\}$.
Since $H$ has minimum degree more than $0.501 t^{1/2}$, each of $S$ and $T$ have size
at least $(0.501 t^{1/2} - 1)^2 > |X_0|/2$, provided $t$ is large enough.
Then since $S,T \subset X_0$ there is a vertex $w \in S \cap T$.
Thus we have paths $xaw$ and $ybw$, where $a \ne b$ since our paths avoid $z$.
However, $xawby$ forms a $5$-cycle, so we have a contradiction.
We conclude $G_t[X_0]$ is empty, and similarly, $G_t[X_1]$ is empty, so $G_t$ is bipartite.

To complete the proof, recall that $e(G_t) \leq (q_t + 1)t/2$ by Proposition \ref{even=}. However by (\ref{tbound}),
$e(G) > (q_t + 1)t/2$ for $t < n$. Thus we must have $t = n$ and $e(G_t) = e(G) = (q + 1)n/2$, so $G_t = G$.
Therefore $G$ itself is bipartite, and by Proposition \ref{even=}, $G$ is the bipartite incidence graph of a projective plane.
\qed

\section{Proof of Theorem \ref{chromatic}}

Let $\ell \geq 2$ and $k \geq 4\ell + 1$ be odd and $c > 0$. Suppose that $G$ is a $\mathcal{C}_{2\ell} \cup \{C_k\}$-free
graph on $n$ vertices with minimum degree at least $cn^{1/\ell}$. We need to show that $\chi(G) < (4k)^{\ell+1}/c^{\ell}$. We use the approach of Thomassen~\cite{Th}. Consider a maximal sequence of vertices $v_1,v_2,\dots,v_s$ such that the $\ell$th neighborhoods $N_{\ell}(v_i)$ are pairwise disjoint. For $r \ge 0$ write $N_{\leq r}(v) = N_0(v) \cup N_1(v) \cup \dots \cup N_{r}(v)$. By Lemma \ref{find-odd}, since $k \ge 2\ell+1$, for any $v \in V(G)$ and $r \leq \ell$ we have $e(N_r(v)) \leq k|N_r(v)|$. Since $G$ is $\mathcal{C}_{2\ell}$-free, no vertex of $N_r(v)$ has more than one neighbor in $N_{r-1}(v)$, so
\[ e(N_{\leq r}(v)) \leq (k + 1)|N_{\leq r}(v)| < 2k|N_{\leq r}(v)|.\]
On the other hand, since $G$ has minimum degree at least $cn^{1/\ell}$,
\[ e(N_{\leq r}(v)) \geq \frac{1}{2}cn^{1/\ell}|N_{\leq r - 1}(v)|.\]
We conclude that $|N_{\le r}(v)| > \frac{1}{4k} cn^{1/\ell}|N_{\le r-1}(v)|$ for all $r \leq \ell$, which implies
\[ |N_{\le \ell}(v)| > \frac{c^{\ell}}{(4k)^{\ell}}n\]
for every vertex $v \in V(G)$. In particular, this holds for $v_1,v_2,\dots,v_s$, so $s < (4k)^{\ell}/c^{\ell}$.
The maximality of $s$ implies that any vertex $v$ is within distance $2\ell$ of some $v_i$. Now we use the assumption that $k \ge 4\ell+1$. The proof of Lemma \ref{find-odd} shows that for $r \le 2\ell$, every subset of $N_r(v)$ induces a subgraph of average degree less than $2k$. It follows that $G[N_r(v)]$ has chromatic number at most $2k$. Furthermore, $G[N_1(v) \cup N_3(v) \cup \dots \cup N_{2\ell - 1}(v)]$ has chromatic number at most $2k$, since there are no edges
between $N_i(v)$ and $N_{i + 2}(v)$ for any $i$, and similarly
$G[N_0(v) \cup N_2(v) \cup \dots \cup N_{2\ell}(v)]$ also has chromatic number at most $2k$. Therefore $G[N_{\leq 2\ell}(v)]$
has chromatic number at most $4k$, so we can cover $N_{\leq 2\ell}(v_i)$ by at most $4k$ independent sets for $1 \le i \le s$. It follows that
\[ \chi(G) \leq 4ks < \frac{(4k)^{\ell + 1}}{c^{\ell}}.\]
This completes the proof of Theorem \ref{chromatic}. \qed

\section{Concluding remarks}

\noindent $\bullet$
Our stability approach not only gives extremal results but describes the approximate structure
of nearly extremal graphs. We only needed these results in the bipartite (Zarankiewicz) setting,
but we note that very similar arguments give analogous results in the non-bipartite (Tur\'an) setting.
For example, we have the following result.
Suppose $G$ is a $C_4$-free graph on $n$ vertices with average degree $d \sim \sqrt{n}$.
Then for any $S, T \sub V(G)$ we have $e(S,T) = \frac{d}{n}|S||T| + o(n^{3/2})$.
The proof is very similar to that of Theorem \ref{bippr}.
First we control the maximum degree as $\Delta<(1+\eps)d$ by deleting $O(n)=o(n^{3/2})$ edges.
Then the argument of Lemma \ref{closedwalks} shows that $w^\circ_6(G) < (1+o(1))n\Delta^2$;
the only difference is that there are $n-1$ choices for $u$ rather than $n/2+o(n)$.
On the other hand, $w^\circ_6(G) = \frac{1}{n}\sum \lambda_i^6$ has a contribution of
$d^6/n \sim n^2$ from the first eigenvalue, so the other eigenvalues are $o(d)$ as $\eps\to 0$.
The pseudorandomness property now follows from the non-bipartite version of Lemma
\ref{pseudo-near-reg-bip}, which is given in \cite[Section 2.4]{KS}.

\noindent $\bullet$  No result similar to Lemma \ref{c4maxdeg} can hold for $C_6$-free graphs:
in fact by the results of \cite{FNV}, there exist $\delta,\eps > 0$ such that any extremal
$C_6$-free graph $G$ with average degree $d$ has at least $\delta e(G)$ edges
 containing a vertex of degree more than $(1 + \eps)d$.

\noindent $\bullet$ We proved the even girth result Theorem \ref{even}
under the assumption that the forbidden odd cycle length satisfies $k \ge 2\ell+3$.
The polarity graphs show that no such result holds for $3 \le k \le \ell$,
but some values of $k$ remain open. In particular, one might think that the
case $k=2\ell+1$ should be approachable by the methods we used to handle $\{C_4,C_5\}$-free
graphs. However, the vertex deletion method does not give sufficient minimum degree
for a straightforward adaptation of this argument, so other ideas are needed.

\noindent $\bullet$ The polarity graphs have arbitrarily large chromatic numbers, as shown by estimates on their independence numbers by Godsil and Newman \cite{GN}. Thus Theorem \ref{chromatic} does not hold when $(k,\ell) \in \{(3,2),(3,3),(3,5),(5,5)\}$,  but it seems likely that it should hold when $k \geq 2\ell + 1$. We also remark that the bound $(4k)^{\ell + 1}/c^{\ell}$ is unlikely to be the correct dependence of $\chi(G)$ on $c$, and it would be interesting to determine this even for $\ell=2$.

\noindent $\bullet$
We remarked earlier that there is little known about the approximate structure
of nearly extremal graphs. In the case of $4$-cycles,
F\"uredi \cite{F1a,F1b} showed that any $C_4$-free graph on $q^2+q+1$ vertices has at most
$\frac{1}{2}q(q+1)^2$ edges, and for large $q$ equality can only hold for polarity graphs.
However, we do not know whether all nearly extremal graphs are structurally close
to polarity graphs. Such an understanding would probably have implications for the conjecture of
Erd\H{o}s that $\mbox{ex}(n,\{C_3,C_4\}) \sim \mbox{z}(n,C_4)$, and for other seemingly more basic questions,
such as whether a graph on $n$ vertices with at least $\mbox{ex}(n,C_4)+1$ edges
must contain many $4$-cycles.

\noindent $\bullet$ We also considered more generally the problem
of determining bipartite graphs $F$ and odd integers $k$ for which $\mbox{ex}(n,\{F,C_k\}) \sim \mbox{z}(n,F)$ and
developed different methods to attack this problem when $F$ is not an even cycle.
This will be the subject of a second paper, co-authored with P. Allen.

\end{document}